\documentclass[a4paper, 11pt]{article}
\usepackage{amsmath} \usepackage{amsfonts} \usepackage{amssymb}\usepackage{latexsym}
\usepackage[english]{babel}
\usepackage{amscd}
\usepackage[dvips,final]{graphics}
\usepackage{locengli,math}
\usepackage{graphicx,pst-plot,psfig,pstricks,pst-node,
pst-text,pst-tree,pst-char,pst-coil,Addaletree}

\begin{document}
\begin{center}
\textbf{\LARGE{\textsf{A simple symmetry generating operads related to rooted planar $m$-ary trees and polygonal numbers}}}
\footnote{
{\it{2000 Mathematics Subject Classification: $05C99,\ 16W30,\ 17A30,\ 17A50, \ 17A42,\ 18D50, \ 55P99$.}}
{\it{Key words and phrases: Free $m$-dendriform algebras, involutive $\mathcal{P}$-Hopf algebras, rooted planar $m$-ary trees, homogeneous polynomials, duality, operads, $k$-gonal numbers.}} 
}
\vskip2cm
\large{
Philippe {\sc Leroux}}
\vskip1cm
To J.-L. Loday for his $60^{th}$ birthday.
\end{center}
\vskip2cm
\noindent
{\bf Abstract:} 
The aim of this paper is to further explore an idea 
from J.-L. Loday briefly exposed in \cite{Lodayscd} and to extend some results obtained in \cite{Loday}. We impose a natural and simple symmetry on a unit action over the most general quadratic relation which can be written. This leads us to two families of binary, quadratic and regular operads whose free objects are computed, as much as possible, as well as their duals in the sense of Ginzburg and Kapranov. Roughly speaking, free objects found here are in relation to $m$-ary trees, triangular numbers and more generally $m$-tetrahedral numbers, homogeneous polynomials on $m$ commutative indeterminates over a field $K$ and polygonal numbers. Involutive connected $\mathcal{P}$-Hopf algebras are constructed and a link to genomics is discussed.
We also propose in conclusion some open questions.

\section{A brief recall on $K$-linear regular operads}
Let $K$ be a null characteristic field. In the sequel, if $S$ is a set, then $KS$ or $K[S]$ will be the free $K$-vector space spanned by $S$. The symbol $S_n$ will denote the symmetric group over
$1, \ldots,n$.
Given a $K$-algebra $A$ `of type $\mathcal{P}$', one considers the family of the $K$-vector spaces $\mathcal{P}(n)$ of $n$-ary operations, see \cite{Lodayren} for instance. We get a linear map
$ \Phi: \mathcal{P}(n) \otimes A^{\otimes n} \xrightarrow{} A, \ \ \Phi(f; (a_1, \ldots, a_n)) \mapsto f(a_1, \ldots, a_n).$
Operations can be composed in the following natural ways. For $f \in \mathcal{P}(m)$, $g \in \mathcal{P}(n)$; $\forall \ 1 \leq i \leq m,$ $f \circ_i g \in \mathcal{P}(m+n-1)$ is defined by:
$$ f \circ_i g (a_1, \ldots, a_{m+n-1}):= f(a_1, \ldots,a_{i-1}, g(a_i, \ldots,a_{i+n-1}), a_{i+n}, \ldots a_{m+n-1}). $$
These composition operations have to obey natural conditions \cite{Lodayren}.
If $h \in \mathcal{P}(l)$, $f \in \mathcal{P}(m)$ and $g \in \mathcal{P}(n)$, then
$ (h \circ_i f) \circ_{j+m-1} g = (h \circ_j g) \circ_i f; \ \ 1\leq i < j\leq l, \ \
 (h \circ_i f) \circ_{i+j-1 } g = h \circ_i (f \circ_j g); \ \ 1 \leq i \leq l; \ 1 \leq j \leq m.$
A $K$-\textit{linear regular operad} $\mathcal{P}$ is then a family of $K$-vector spaces $(\mathcal{P}(n))_{n>0}$ equipped with composition maps $\circ_i$ verifying the above relations. If all possible operations are generated by composition from $\mathcal{P}(2)$, then the operad is said to be {\it{binary}}. It is said to be {\it{quadratic}} if all the relations between operations are consequences of relations described exclusively with the help of monomials with two operations. Ideals of such an algebra $A$ of type $\mathcal{P}$
is a submodule of it such that for any $\diamond \in \mathcal{P}(2)$, $x \diamond y \in I$ whenever one of the variables is in $I$. Clearly, the quotient $A/I$ is still an algebra $A$ of type $\mathcal{P}$.
Let $V$ be a $K$-vector space. The {\it{free $\mathcal{P}$-algebra}} $\mathcal{P}(V)$ on $V$ is by definition
a $\mathcal{P}$-algebra equipped with a linear map $i: \ V \rightarrow \mathcal{P}(V)$ which satisfies the following universal property:
for any linear map $f: V \xrightarrow{} A$, where $A$ is a $\mathcal{P}$-algebra, there exists a unique  $\mathcal{P}$-algebra morphism $\phi: \mathcal{P}(V) \xrightarrow{} A$ such that $\phi \circ i=f$.
Since our $\mathcal{P}$-algebras are regular, the free $\mathcal{P}$-algebra over a $K-$vector space $V$ is of the form:
$ \mathcal{P}(V) := \bigoplus_{n \geq 1} \ \mathcal{P}_n \otimes V^{\otimes n} ,$ with $\mathcal{P}(n)=\mathcal{P}_n \otimes K[S_n]$.
In particular, the free $\mathcal{P}$-algebra on one generator is $ \mathcal{P}(K) := \bigoplus_{n \geq 1} \mathcal{P}_n$.
Recall that in the case of a regular binary and quadratic operad $\mathcal{P}$, the free $\mathcal{P}$-algebra is entirely induced by the free $\mathcal{P}$-algebra on one generator. The generating function of the regular operad $\mathcal{P}$, or its Poincar\'e series, is given by:
$ f^{\mathcal{P}}(x):= \sum \ (-1)^n \textrm{dim} \ \mathcal{P}_n \ x^n.$
Below, we will indicate the sequence $(\textrm{dim} \ \mathcal{P}_n)_{n \geq 1}$.
Let $\mathcal{P}$ be a binary quadratic operad. By a \textit{unit action} \cite{Lodayscd}, we mean the choice of two linear applications:
$\upsilon: \mathcal{P}(2) \xrightarrow{} \mathcal{P}(1)$ and $\varpi:\mathcal{P}(2) \xrightarrow{} \mathcal{P}(1),$
giving sense, when possible, to $x \diamond 1$ and $1 \diamond x$, for all operations $\diamond \in \mathcal{P}(2)$ and for all $x$ in the $\mathcal{P}$-algebra $A$, i.e.,
$x \diamond 1 = \upsilon(\diamond)(x)$ and $1 \diamond x= \varpi(\diamond)(x)$.
If $\mathcal{P}(2)$ contains an associative operation, say $\star$, then we require that $x \star 1 := x := 1 \star x$, i.e., $\upsilon(\star) := Id := \varpi(\star)$.
We say that the unit action, or the couple $(\upsilon,\varpi)$ is {\it{compatible}} with the relations of the $\mathcal{P}$-algebra $A$ if they still hold on $A_+:= K.1 \oplus A$ as far as the terms are defined.
Let $A$, $B$ be two $\mathcal{P}$-algebras such that $\mathcal{P}(2)$ contains an associative operation $\star$. Using the couple $(\upsilon,\varpi)$, we extend binary operations $\diamond \in \mathcal{P}(2)$ to the $K$-vector space $A \otimes 1.K \oplus K.1 \otimes B \oplus A \otimes B$ by requiring:
\begin{eqnarray}
(a \otimes b) \diamond (a' \otimes b') &:= &(a \star a') \otimes (b \diamond b') \ \ \ \textrm{if} \ \ b \otimes b' \not= 1 \otimes 1, \\
(a \otimes 1) \diamond (a' \otimes 1) &:= &(a \diamond a') \otimes 1, \ \ \ \textrm{otherwise}.
\end{eqnarray}
The unit action or the couple $(\upsilon,\varpi)$ is said to be {\it{coherent}} with the relations of $\mathcal{P}$ if $A \otimes 1.K \oplus K.1 \otimes B \oplus A \otimes B$, equipped with these operations is still a $\mathcal{P}$-algebra. Observe that a necessary condition for having coherence is compatibility.
\begin{theo} [\textbf{Loday \cite{Lodayscd}}]
\label{Lodaych}
Let $\mathcal{P}$ be a binary quadratic operad. Suppose there exists an associative operation in $\mathcal{P}(2)$. Then, any unit action coherent with the relations of $\mathcal{P}$ equips the augmented free $\mathcal{P}$-algebra $\mathcal{P}(V)_+$ on a $K$-vector space $V$ with a
coassociative coproduct $\Delta: \mathcal{P}(V)_+ \xrightarrow{}  \mathcal{P}(V)_+ \otimes \mathcal{P}(V)_+,$
which is a $\mathcal{P}$-algebra morphism.
Moreover,
$\mathcal{P}(V)_+$ is a connected Hopf algebra.
\end{theo}
Such structures are often called connected $\mathcal{P}$-Hopf algebras to indicate that the coproduct is a morphism for any operations in $\mathcal{P}(2)$.
\section{On consequences of a simple symmetry on a unit action}
\noindent
\textbf{Motivation:} We further explore an idea briefly exposed in \cite{Lodayscd}.
Let $k>1$ be an integer. We consider $k$ binary operations $\bullet_i$, $1\leq i \leq k$, over a $K$-vector space $V$, supposed to be related by quadratic and regular relations. We suppose the existence of a unit element, denoted by 1, which acts as follows. First we rename $\bullet_1$ as $\succ$ and $\bullet_n$ as $\prec$. Second, set for all $x \in V$, $1 \prec x=0=x\succ 1$, $x\prec 1=x= 1 \succ x$ and for all $2 \leq i \leq k-1$ set $x \bullet_i 1=0= 1 \bullet_i x$. Observe that the action of the unit is invariant under the transformation,
\begin{eqnarray}
\label{transf} 
x\prec y \longmapsto y \succ x,\hskip1cm 
x\succ y \longmapsto y \prec x,\hskip1cm 
x\bullet_i y \longmapsto y \bullet_i x,
\end{eqnarray} 
for all $2 \leq i \leq k-1$.
We now write down the most general equation relating operations $\bullet_i$, $1 \leq i \leq k$  to their braces, that is:
$$ \sum_{i,j=1}^k \lambda_{ij} \ (x \bullet_i y) \bullet_j z=\sum_{i,j=1}^k \lambda'_{ij} \ x \bullet_i (y \bullet_j z), $$
where the $\lambda_{ij}$ and  $\lambda'_{ij}$ are scalars of $K$. Setting $x=1$, then $y=1$, then $z=1$ in the previous equation and applying our choice of the unit action lead to the following system of relations. For all $x,y,z \in V,$ for all $2 \leq i \leq k-1$, we get:
$$ (x\prec y)\prec z = x \prec (y \star z), \ \  \  \  \ (x \prec y)\bullet_i z =x\bullet_i (y \succ z) $$ 
$$ (x \succ y)\prec z = x \succ (y \prec z), \ \  \  \  \ (x \succ y)\bullet_i z =x\succ (y \bullet_i z) $$
$$(x \star y)\succ z =x\succ (y \succ z) , \ \  \  \  \ (x\bullet_i  y)\prec z = x \bullet_i (y \prec z), $$
where $x \star y := x\prec y + x\succ y$. Observe that $\star$ is associative and that this system of $3(k-1)$ relations is invariant under the transformation (\ref{transf}). 
Binary, regular and quadratic operads denoted by $^k\mathcal{P}$ can be naturally associated with each system of relations. The coefficients of their Poincar\'e series starts with $1,k, 2k^2-3(k-1), \ldots$, that is for
$k=2$, we get $1,2,5, \ldots$ which is the beginning of the Catalan numbers counting for instance the number of planar rooted binary trees on $p$ internal vertices, with $p>0$. The operad $^2\mathcal{P}$ is the operad $Dend$ of dendriform algebras introduced by J.-L. Loday \cite{Loday}.
For $k=3$, we get $1,3,12, \ldots$ which is the beginning of the sequence counting the number of planar rooted ternary trees on $p$ internal vertices, with $p>0$. This operad will be investigated bellow. For $k=4$, we get $1,4,23, \ldots$ which is the beginning of the sequence counting the number of non-crossing 
connected graphs \cite{Flajolet} on $p+1$ vertices, with $p>0$ and so on. 
By the existence of an involution on $^k\mathcal{P}(V)$, the free $^k\mathcal{P}$-algebra over $V$, we mean a linear involutive map $\dagger: V \longrightarrow V$  acting as follows,
\begin{eqnarray} 
\label{inv}
(x\prec y)^{\dagger}= y^{\dagger}\succ x^{\dagger}, \ \ (x\succ y)^{\dagger}= y^{\dagger}\prec x^{\dagger}, \ \ (x\bullet_i y)^{\dagger}= y^{\dagger}\bullet_{n+1-i} x^{\dagger},
\end{eqnarray}
for $2 \leq i \leq k-1$.
Extend $\dagger$ on $^k\mathcal{P}(V)^{\otimes 2}$
by the following formula,
$ (x \otimes y) ^{\dagger}= x^{\dagger} \otimes y^{\dagger}.$
An involutive connected $\mathcal{P}$-Hopf algebra is a connected $\mathcal{P}$-Hopf algebra equipped with an involution $\dagger$ such that $\Delta(x^\dagger)=\Delta(x)^\dagger$. 
Then, following theorems are direct applications of Theorem \ref{Lodaych}.
\begin{theo}
\label{hopf}
Let $V$ be a $K$-vector space.
For all $n>1$, there exists a connected $\mathcal{P}$-Hopf algebra structure on $^k\mathcal{P}(V)$, the free $^k\mathcal{P}$-algebra over $V$, which is involutive if $^k\mathcal{P}(V)$ has an involution.
\end{theo}
\Rk \textbf{(Opposite  $^k\mathcal{P}$-algebra.)}
Let $T$ be a $^k\mathcal{P}$-algebra. Define new operations by:
$$x \prec' y:= y \succ x; \ x \succ' y:= y \prec x; \ x \bullet_i' y:= y \bullet_{n+1-i} x,$$
for all $2 \leq i \leq k-1$. Then, the $K$-vector space $T$ equipped with these operations is a new
$^k\mathcal{P}$-algebra, denoted by $T^{op}$, called the opposite $^k\mathcal{P}$-algebra. A $^k\mathcal{P}$-algebra is said to be commutative if $T^{op}=T$. 
\begin{theo}
Let $V$ be a $K$-vector space.
For all $k>1$, there exists a connected $\mathcal{P}$-Hopf algebra structure on $^k\mathcal{P}_{com}(V)$, the free commutative $^k\mathcal{P}$-algebra over $V$, which is involutive if $^k\mathcal{P}_{com}(V)$ has an involution.
\end{theo}
\Proof
Observe that the choice of the unit action is such that $x \prec' 1:= 1 \succ x:=x$ and $x \succ' 1:= 1 \prec x:=x$ and thus in agreement with the opposite structure.
\eproof

The main aim of this paper is to focus on systems of 
quadratic relations exposed just above and giving birth to regular, binary and quadratic operads, to find
the free object associated with each operad, to compute the associated dual operad in the sense of Ginzburg and Kapranov \cite{GK} and the free object
associated with it as much as possible. We start with Section~3 by exploring the case over three operations. Roughly speaking, the associated operad
is called $3-Dend.$ and the free object can be constructed over planar rooted ternary trees. The operadic dual is explored in Section~4 and is called $Triang.$ since its structure is related to triangular numbers. It is shown that the free triangular algebra on one generator can be constructed over homogeneous polynomials in three commutative indeterminates, so is in relation to projective algebraic curves in the projective plane $\mathbb{P}^2(K)$. In both cases,
a chain-complex is constructed and an explicit notion of (co)homology is proposed.    
To extend these results to planar rooted $m$-ary trees, we add to our systems of quadratic relations 
other natural quadratic constrains. This is explored in Section~5. The operadic dual (in the sense of Ginzburg and Kapranov) is computed, giving birth to the operad $m-Tetra.$ 
whose structure is related to tetrahedral numbers in dimension $m-1$. Here again, it is shown that the free $m$-tetrahedral algebra on one generator can be constructed over homogeneous polynomials in $m$ commutative indeterminates, and so is related to projective algebraic hypersurfaces
in $\mathbb{P}^{m-1}(K)$. A (Co)homology theory for these categories is also proposed and a link to genomics is also discussed. We gather in Section~6 all the coefficients of the Poincar\'e series of the operads $m-Tetra.$ in a triangle shape giving birth
to the well-known Pascal triangle. We propose then a dual triangle of the Pascal's one by gathering this time the coefficients of the Poincar\'e series of the dual operads $m-Dend.$.
After triangular shapes and $m$-ary trees, we compute in Section~7,
the operadic dual of systems of quadratic relations
exposed just above and prove that free objects associated with them are related to $k$-gonal numbers. We do not know for the time being how to construct systematically the free objects associated with $^k \mathcal{P}$, $k \geq 4$.
\section{An operad over rooted planar ternary trees }
We start with the case of rooted planar ternary trees to expose our results which will
be generalised in the other sections. We voluntarily adopt the strategy of J.-L. Loday \cite{Loday} to tackle the issue.
\subsection{On rooted planar ternary trees}
We now focus on the case $n=3$ to show that the free object associated with the operad $^3\mathcal{P}$
can be explicited thanks to the planar rooted ternary trees (ternary trees for short). For convenience, we rename $^2\mathcal{P}$, that is $Dend$ in $2-Dend$
and $^3\mathcal{P}$ in $3-Dend$.
\begin{defi}{}
\label{def3}
A $K$-vector space $T$ is a 3-dendriform algebra if it is equipped with 3 binary operations $\prec,\succ,\cdot: T^{\otimes 2} \longrightarrow T$ verifying for all $x,y,z \in T,$
$$ 1. \ (x\prec y)\prec z = x \prec (y \star z), \ \  \  \  \ 4. \ (x \prec y)\bullet z =x\bullet (y \succ z) $$ 
$$ 2. \ (x \succ y)\prec z = x \succ (y \prec z), \ \  \  \   5. \ (x \succ y)\bullet z =x\succ (y \bullet z) $$
$$3. \ (x \star y)\succ z =x\succ (y \succ z) , \ \  \  \  \ 6. \ (x\bullet  y)\prec z = x \bullet (y \prec z), $$
where $x \star y := x\prec y + x\succ y$. It is said to be involutive if there exists a linear involutive map
$\dagger:T \rightarrow T$ verifying $(x\prec y)^\dagger=y^\dagger \succ x^\dagger$, $(x\succ y)^\dagger=y^\dagger \prec x^\dagger$ and $(x\bullet y)^\dagger=y^\dagger \bullet x^\dagger$, for all $x,y \in T$.
\end{defi}
Morphisms 
of 3-dendriform algebras are straightforward to define. The functorial diagram between categories holds:
\begin{center}
$
\begin{array}{ccc} 
2-\textsf{Dend.} &  \longrightarrow  & 3-\textsf{Dend.}  \\ 
_+ \searrow&       &\swarrow _+  \\
&  \textsf{As.}   & \\ 
\end{array} 
$
\end{center}
For more information about ternary trees, we refer to the extended literature. In small dimensions, one gets:
$\treeg_0:=\{\mid \}$,$\treeg_1:=\{\treeg \}$, 
$\treeg_2:=\{\treegg, \ \treegm \ \treedd \}$,
$$\treeg_3:=\{\treeggg, \
\treegmm, \
\treeggd, \
\treeggm, \
\treegmd, \
\treemgm, \
\treeddd, \
\treeGm, \
\treegdg, \
\treedmm, \
\treemmd, \
\treemmg\}$$
and so on, where $\treeg_n$
denotes the set of ternary trees of degree $n$ , that is with $n$ internal vertices and thus with $2n+1$
leaves. The cardinal of $\treeg_n$ is equal to $\frac{(3n)!}{n!(2n+1))!}$.
Each tree of $\treeg_n$, $n>0$, can be decomposed in a unique way \textit{via} the so called grafting operation (which is ternary). It is defined for all $p,q,r \in \mathbb{N}$ as follows, 
$$\vee:\treeg_p \times \treeg_q \times \treeg_r \longrightarrow \treeg_{p+q+r+1},$$
$$ (t_1,t_2,t_3) \longmapsto t_1\vee t_2 \vee t_3,$$
the symbol $t_1\vee t_2 \vee t_3$ meaning that the root of $t_1$ is glued to the left leave of $\treeg$, the root of $t_2$ is glued to its middle leave and the root of $t_3$ to its right leave. For instance, $\treegm= \mid \vee \ \treeg \ \vee \mid$.

There exists also an involution denoted by $\dagger:\treeg_n \mapsto \treeg_n$ for all $n$, defined recursively by $t^\dagger=t_3^\dagger\vee t_2^\dagger \vee t_1^\dagger$ if $t=t_1\vee t_2 \vee t_3$. Pictorially, this is the mirror image \textit{via} its central axis. For instance, $\treegg \ ^\dagger =\treedd$.

Over the $K$-vector space $K[\treeg_\infty]:=K\mid \oplus \hat{K[\treeg_\infty]}$, where $\hat{K[\treeg_\infty]}:=\bigoplus_{n>0} K\treeg_n$, we introduce recursively the following binary operations first on the trees. They are naturally extended by bilinearity to the whole $\hat{K[\treeg_\infty]}$. Let $p,q>0$ and
set for any $t=t_1\vee t_2 \vee t_3 \in \treeg_p$ and $r=r_1\vee r_2 \vee r_3 \in \treeg_q$,
\begin{eqnarray}
t\succ r &=& (t\star r_1)\vee r_2 \vee r_3, \\ 
t\prec r &=& t_1\vee t_2 \vee (t_3\star r),\\
t\bullet r &=& r_1\vee (r_2 \vee (r_3 \star t_1)\vee t_2) \vee t_3,\\
\mid \prec t &=& 0= t\succ \mid,\hskip1cm
t \prec \mid = t= \mid \succ t,\\
\mid\bullet t &=& 0=t\bullet \mid,
\end{eqnarray} 
where as usual $x \star y= x\prec y + x \succ y$ is
by construction associative. The symbols $1\prec 1$, $1\succ 1$ and $1 \bullet 1$ are not defined. Observe that $
\mid \star t = t= t\star \mid$ and that our three operations repect the natural grading of $\hat{K[\treeg_\infty]}$ since $\bullet,\succ,\prec: K\treeg_p \otimes K\treeg_q \longrightarrow K\treeg_{p+q}$.
\begin{prop}
Equipped with these three binary operations, $\hat{K[\treeg_\infty]}$ is an involutive 3-dendriform algebra generated by $\treeg$.
\end{prop}
\Proof
We proceed by induction on the degree of trees.
Observe that the involution $\dagger$ on $\treeg_2$
acts as expected since $(\treeg \prec \treeg)^\dagger=\treeg \succ \treeg $ and so on.
By induction, one checks that $\dagger$ is an involutive map see \ref{def3}. For instance, if $r=r_1 \vee r_2 \vee r_3$ and $t$ are trees, then
$(r \prec t)^\dagger = (r_1 \vee r_2 \vee (r_3 \star t))^\dagger=(r_3 \star t)^\dagger \vee r_2^\dagger \vee r_1^\dagger$ by definition and
$(r_3 \star t)^\dagger \vee r_2^\dagger \vee r_1^\dagger=(t^\dagger \star r_3^\dagger) \vee r_2^\dagger \vee r_1^\dagger $ by induction. Therefore, $(r \prec t)^\dagger =  r^\dagger \succ t^\dagger$. Let us check the 6 axioms of Definition \ref{def3}. Let $r,s,t$ be ternary trees. We get:
$(r \prec s) \prec t = r_1 \vee r_2 \vee (r_3\star s)\star t=
r_1 \vee r_2 \vee r_3 \star (s\star t)$ by induction, therefore Axiom 1 holds and Axiom 3 as well \textit{via} involution on Axiom 1. Axiom 2 is straightforward.
Axiom 4 leads to:
$(r \prec s) \bullet t= r_1 \vee r_2 \vee (r_3\star s) \bullet t= r_1 \vee (r_2 \vee ((r_3\star s) \star t_1) \vee t_2) \vee t_3= r_1 \vee (r_2 \vee (r_3\star (s \star t_1)) \vee t_2) \vee t_3$ by induction, hence Axiom 4 holds. Axiom 5 and 6 are straightforward. Therefore, $\hat{K[\treeg_\infty]}$ is an involutive 3-dendriform algebra. We introduce the middle map $m: \treeg_p \longrightarrow \treeg_{p+1}$ for all $p \in \mathbb{N}$, such that $t \mapsto \mid \vee t \vee \mid$. We now prove that $\hat{K[\treeg_\infty]}$ is generated by $\treeg$ by induction on the degree of trees.
Indeed, the result holds in small dimension (checked by hand up to $p=3$). Moreover we have for a tree $t$,
\begin{eqnarray*}
t:=t_1 \vee t_2 \vee t_3 &=& \treeg \ \ \ \ \textrm{if} \ \ t_1=t_2=t_3=\mid, \\
&=&  \treeg \prec t_3 \  \ \textrm{if} \ \ t_1=t_2=\mid, \ t_3 \not=\mid, \\
&=& t_1 \succ \treeg  \ \ \textrm{if} \ \  t_1\not=\mid,\ t_2=\mid= t_3 , \\
&=&t_1\succ (m((t_2)_1) \bullet ((t_2)_2 \succ m((t_2)_3)))\prec t_3  \ \ \textrm{otherwise}. 
\end{eqnarray*} 
\eproof
\begin{prop}
\label{x}
The unique 3-dendriform algebra map $3-Dend(K) \longrightarrow \hat{K[\treeg_\infty]}:= \bigoplus_{n>0} K[\treeg_n]$ sending the generator $x$ of $3-Dend(K)$ to $\treeg$ is an isomorphism. 
\end{prop}
\Proof
We will check that $(\hat{K[\treeg_\infty]},\prec,\succ,\bullet)$ verifies the universal condition to be the free 3-dendriform algebra on one generator. Let $T$ be a 
3-dendriform algebra and let $a\in T$. By induction, we construct a linear map $\alpha:\hat{K[\treeg_\infty]} \longrightarrow T$ on its values on ternary trees as follows. Let $t=t_1 \vee t_2 \vee t_3 \in \treeg_p$ and set:
\begin{eqnarray*}
\alpha(t_1 \vee t_2 \vee t_3) &=& a \ \ \ \ \textrm{if} \ \ t_1=t_2=t_3=\mid, \\
&=&  a \prec \alpha(t_3) \  \ \textrm{if} \ \ t_1=t_2=\mid, \ t_3 \not=\mid, \\
&=& \alpha(t_1) \succ a  \ \ \textrm{if} \ \  t_1\not=\mid,\ t_2=\mid= t_3 , \\
&=&\alpha(t_1)\succ (\alpha(m((t_2)_1)) \bullet (\alpha((t_2)_2) \succ \alpha(m((t_2)_3))))\prec \alpha(t_3)  \ \ \textrm{otherwise}. 
\end{eqnarray*} 
The map $\alpha$ is unique since  $\hat{K[\treeg_\infty]}$ is generated by $\treeg$ and that $\alpha(\treeg)=a$. It is a morphism of 3-dendriform algebras as one can show by induction on the degree of trees. For instance,
\begin{eqnarray*}
\alpha(r\prec t)&:=&\alpha(r_1)\succ (\alpha(m((r_2)_1)) \bullet (\alpha((r_2)_2) \succ \alpha(m((r_2)_3))))\prec \alpha(r_3 \star t)\\
 &=& \alpha(r_1)\succ (\alpha(m((r_2)_1)) \bullet (\alpha((r_2)_2) \succ \alpha(m((r_2)_3))))\prec (\alpha(r_3) \star \alpha(t)) \ \ \textrm{by induction}\\
&=&(\alpha(r_1)\succ (\alpha(m((r_2)_1)) \bullet (\alpha((r_2)_2) \succ \alpha(m((r_2)_3))))\prec \alpha(r_3)) \prec \alpha(t) \ \ \textrm{by Axiom 1}\\
&:=& \alpha(r) \prec \alpha(t).
\end{eqnarray*}
The relation $\alpha(r \bullet t)=\alpha(r) \bullet \alpha(t)$ follows from the following general equality: Let
$a,b,c,a',b',c' \in T$, where $T$ is any 3-dendriform algebra. Then,
\begin{eqnarray*}
(a\succ b \prec c)\bullet (a'\succ b' \prec c')=
a\succ(b\bullet((c\star a')\succ b')\prec c',
\end{eqnarray*}
holds.
Therefore, $(\hat{K[\treeg_\infty]},\prec,\succ,\bullet)$ is the free dendriform algebra on one generator. 
\eproof
\begin{theo}[Free 3-dendriform algebra]
Let $V$ be a $K$-vector space.
The unique 3-dendriform algebra map $3-Dend(V) \longrightarrow \bigoplus_{n>0} K[\treeg_n]\otimes V^{\otimes n},$ which sends the generator $v\in V$ to $\treeg \otimes v$ is an isomorphism.
\end{theo}
\Proof
Define on $\bigoplus_{n>0} K[\treeg_n]\otimes V^{\otimes n}$ the following 3-dendriform algebra structure:
\begin{eqnarray*}
t\otimes \omega \prec t'\otimes \omega' &:=& t\prec t'\otimes \omega\omega', \\
t\otimes \omega \succ t'\otimes \omega' &:=& t\succ t'\otimes \omega\omega', \\
t\otimes \omega \bullet t'\otimes \omega' &:=& t\bullet t'\otimes \omega\omega'. \\
\end{eqnarray*}
Since the relations defining 3-dendriform algebras
are regular, the free 3-dendriform algebra over $V$ is then determined by the free 3-dendriform algebra on one generator, hence $3-Dend(V):=\bigoplus_{n>0} 3-Dend(K)_n \otimes V^{\otimes n}=\bigoplus_{n>0} K[\treeg_n]\otimes V^{\otimes n}$ by Proposition \ref{x}.
\eproof
\begin{coro}
The augmented free 3-dendriform algebra on the generator $\treeg$, $K[\treeg_\infty]$, is naturally equipped with an involutive connected $\mathcal{P}$-Hopf algebra.
\end{coro}
\Proof
Even if this result is a corollary of Theorem \ref{hopf}, we detail the proof to help the reader to understand Theorem \ref{Lodaych}.
Extend the binary operations $\prec,\succ,\bullet$ of
$\hat{K[\treeg_\infty]}$ to $K[\treeg_\infty]$ as follows:
$ t \prec 1 := t, \  \ 1 \prec t := 0, \ \ 1 \succ t := t, \ \ t \succ 1 := 0, \ \ 1\bullet t=0=t \bullet 1$ for all $t \in \treeg_n$, $n\not=0$, next by linearity. We cannot extend the operations $\prec$, $\succ$ and $\bullet$ to $K$, i.e., $1 \prec 1$, $ 1 \succ 1$ and $1 \bullet 1$ are not defined.
Let us show that this choice is compatible. Let $x,y,z \in K[\treeg_\infty]$. We have to show for instance that the relation $(x \prec y)\prec z = x\prec (y \star z)$ holds in $K[\treeg_\infty]$. Indeed for $x=1$, we get $0=0$. For $y =1$ we get $ x \prec z = x \prec z$ and for $z=1$ we get $ x \prec y = x \prec y$. We do the same thing with the 5 other equations to find that the augmented free 3-dendriform algebra $K[\treeg_\infty]$ is still a 3-dendriform algebra.
Secondly, let us show that this choice is coherent. Let $x_1,x_2, x_3, y_1, y_2, y_3 \in K[\treeg_\infty]$. We have to show that, for instance:
$(*) \ \ \  ((x_1 \otimes y_1) \prec (x_2 \otimes y_2)) \prec (x_3 \otimes y_3) = (x_1 \otimes y_1) \prec ((x_2 \otimes y_2) \star (x_3 \otimes y_3)).$
Indeed, if there exists a unique $y_i =1$, the others belonging to $K[\treeg_\infty]$, then we get:
$ x_1 \star x_2 \star x_3  \otimes (y_1 \prec y_2) \prec  y_3 = x_1 \star x_2 \star x_3  \otimes y_1 \prec ( y_2 \star  y_3),$
which always holds since our choice of the unit action is compatible. Similarly if $y_1=y_2=y_3=1$.
If $y_1=y_2=1$ and $y_3 \in K[\treeg_\infty]$, we get: $0=0$, similarly if $y_1=1=y_3$ and $y_2 \in K[\treeg_\infty]$. If $y_1 \in K[\treeg_\infty]$ and $y_2=1=y_3$, the two hand sides of $(*)$ are equal to $x_1 \star x_2 \star x_3  \otimes y_1$. Therefore $(*)$ holds in $\hat{K[\treeg_\infty]}  \otimes 1.K \oplus K.1 \otimes \hat{K[\treeg_\infty]}  \oplus \hat{K[\treeg_\infty]}  \otimes \hat{K[\treeg_\infty]} $. Checking  the 5 other relations shows that our choice of the unit action is coherent.
From Theorem \ref{Lodaych}, we recover a connected $\mathcal{P}$-Hopf algebra structure on the augmented free 3-dendriform algebra generated by $\treeg$. 
Extend the involutive map on $\hat{K[\treeg_\infty]}^{\otimes 2}$ by setting
$(x \otimes y)^{\dagger}:=x^{\dagger} \otimes y^{\dagger},$ for any $x,y \in \hat{K[\treeg_\infty]}$. Then, the connected $\mathcal{P}$-Hopf algebra with coproduct $\Delta$ just defined is involutive since $\Delta (x^\dagger) = \Delta(x)^\dagger.$
\eproof
\subsection{(Co)Homology of 3-dendriform algebras}
We show the existence of a chain complex of Hochschild type for any 3-dendriform algebra and define a (co)homology theory for this category. Let $T$ be a 3-dendriform algebra. For all $n \in \mathbb{N}$, define $X_n:=\{(k,j); 0\leq k \leq j \leq n\}$ and the module of $n$-chains of $T$ as 
$C_n ^{3-Dend}(T):=K(X_n)\otimes T^{\otimes n}$. Introduce the differential operator $d:=\sum_{1 \leq i \leq n-1} \ (-1)^{i+1}d_i$, where
$d_i:C_n ^{3-Dend}(T)\longrightarrow C_{n-1} ^{3-Dend}(T)$, $1 \leq i \leq n-1$, are the face operators and act on $X_n$ as follows:

\noindent
\textbf{Step 1:} Set $d_i(k,j):=(\tilde{d_i}(k), \tilde{d_i}(j))$ where
$\tilde{d_i}:\{1,\ldots, n\} \longrightarrow \{1,\ldots, n-1 \}$ is such that $\tilde{d_i}(r):=r-1$ 
if $i \leq r$ and $\tilde{d_i}(r):= r$ if $i \geq r+1$. The face maps are extended linearly to maps
$d_i:K[X_n] \longrightarrow K[X_{n-1}]$.

\noindent
\textbf{Step 2:} Introduce now the symbol:
\begin{eqnarray*}
\circ_i^{(k,j)}:=
\begin{cases}
\bullet \ \ \ \ if \  \ i-1\in \{k,j\}; \ i \in \{k,j\},\\
\succ \ \ if \ \ i-1\notin \{k,j\}; \ i \in \{k,j\},\\
\prec \ \ if \ \ i-1\in \{k,j\}; \ i \notin \{k,j\},\\
\star \ \ \ \ if \ \ i-1 \notin \{k,j\}; \ i \notin \{k,j\}.
\end{cases}
\end{eqnarray*}
We now explicit the action of the face maps $d_i:C_n ^{3-Dend}(T)\longrightarrow C_{n-1} ^{3-Dend}(T)$, $1 \leq i \leq n-1$ by $d_i((k,j); x_1 \otimes \ldots \otimes x_n):= (d_i(k,j);x_1 \otimes \ldots \otimes x_{i-1} \otimes x_i \circ_i^{(k,j)} x_{i+1}\otimes \ldots \otimes x_n)$.
\begin{prop}
The face maps $d_i$, $1 \leq i \leq n-1$ satisfy the simplicial relations $d_id_j=d_{j-1}d_i$ for $i<j$. Therefore $(C_*^{3-Dend}(T),d)$ is a chain-complex.
\end{prop}
\Proof
We prove first that $d_1d_2=d_1d_1$ on $C_3^{3-Dend}(T)$. Let $x,y,z \in T$.
\begin{eqnarray*}
\begin{cases}
d_1d_2((0,0);x \otimes y \otimes z)=d_1((0,0);x \otimes y \star z)=((0,0);x \prec (y \star z)),\\
d_1d_1((0,0);x \otimes y \otimes z)=d_1((0,0);x \prec y \otimes z)=((0,0);(x \prec y) \prec z),
\end{cases}
\end{eqnarray*}
hence the equality \textit{via} Axiom 1.
\begin{eqnarray*}
\begin{cases}
d_1d_2((1,1);x \otimes y \otimes z)=d_1((1,1);x \otimes y \prec z)=((0,0);x \succ (y \prec z)),\\
d_1d_1((1,1);x \otimes y \otimes z)=d_1((0,0);x \succ y \otimes z)=((0,0);(x \succ y) \prec z),
\end{cases}
\end{eqnarray*}
hence the equality \textit{via} Axiom 2.
\begin{eqnarray*}
\begin{cases}
d_1d_2((2,2);x \otimes y \otimes z)=d_1((1,1);x \otimes y \succ z)=((0,0);x \succ (y \succ z)),\\
d_1d_1((2,2);x \otimes y \otimes z)=d_1((1,1);x \star y \otimes z)=((0,0);(x \star y) \succ z),
\end{cases}
\end{eqnarray*}
hence the equality \textit{via} Axiom 3.
\begin{eqnarray*}
\begin{cases}
d_1d_2((0,2);x \otimes y \otimes z)=d_1((0,1);x \otimes y \succ z)=((0,0);x \bullet (y \succ z)),\\
d_1d_1((0,2);x \otimes y \otimes z)=d_1((0,1);x \prec y \otimes z)=((0,0);(x \prec y) \bullet z),
\end{cases}
\end{eqnarray*}
hence the equality \textit{via} Axiom 4.
\begin{eqnarray*}
\begin{cases}
d_1d_2((1,2);x \otimes y \otimes z)=d_1((1,1);x \otimes y \bullet z)=((0,0);x \succ (y \bullet z)),\\
d_1d_1((1,2);x \otimes y \otimes z)=d_1((0,1);x \succ y \otimes z)=((0,0);(x \succ y) \bullet z),
\end{cases}
\end{eqnarray*}
hence the equality \textit{via} Axiom 5.
\begin{eqnarray*}
\begin{cases}
d_1d_2((0,1);x \otimes y \otimes z)=d_1((0,1);x \otimes y \prec z)=((0,0);x \bullet (y \prec z)),\\
d_1d_1((0,1);x \otimes y \otimes z)=d_1((0,0);x \bullet y \otimes z)=((0,0);(x \bullet y) \prec z),
\end{cases}
\end{eqnarray*}
hence the equality \textit{via} Axiom 6.
The sequel of the proof splits into two cases. The case $j>i+1$ is straightforward and the case $j=i+1$ depends on the computations above and the fact that $\star$ is associative.
\eproof
\section{The operad \textit{Triang.}}
To explain how we obtain the (co)homology of 3-dendriform algebras, we need to study the dual operad in the sense of Ginzburg and Kapranov. We will
show that its Poincar\'e series is $f_{Triang.}(x)=\sum_{n=1}^{\infty} \ (-1)^n \frac{n(n+1)}{2}x^n=-\frac{x}{(x+1)^3}$, where the triangular number, $\frac{n(n+1)}{2}$, is the dimension of $Triang._n$ (see introduction), hence the name of this operad.
\begin{defi}{}
A triangular algebra $T$ is a $K$-vector space equipped with three binary operations $\perp,\vdash,\dashv: T^{\otimes 2} \longrightarrow
T$ verifying for all $x,y,z \in T$, the following 12 axioms.

\begin{eqnarray}
\begin{cases}
1. \ (x\dashv y)\dashv z= x\dashv(y\dashv z),& 
6. \ (x\dashv y)\perp z=x \perp(y\vdash z),\\
2. \ (x\dashv y)\dashv z= x\dashv(y\vdash z),& 
7. \ (x\vdash y)\perp z=x \vdash(y\perp z),\\
3. \ (x\vdash y)\dashv z= x\vdash(y\dashv z),& 
8. \ (x\perp y)\dashv z=x \perp (y\dashv z),\\
4. \ (x\dashv y)\vdash z= x\vdash(y\vdash z),& 
9. \ (x\perp y)\perp z=0 \overset{10.}{=}  x \perp(y\perp z),\\
5. \ (x\vdash y)\vdash z= x\vdash(y\vdash z),& 
11. \ (x\perp y)\vdash z=0 \overset{12.}{=} x \dashv(y\perp z).
\end{cases}
\end{eqnarray}
\end{defi}
\begin{exam}{}
Any associative algebra $(A, \cdot)$ is a triangular algebra by setting $\vdash = \cdot =\dashv$ and $\perp=0$. In the non-graded setting,
let $(A,d)$ be a differential associative algebra, that is $d(ab)=d(a)b + ad(b)$ and $d^2=0$. Set
$a\dashv b := ad(b)$, $a\vdash b := d(a)b$ and
$a \perp b:=d(a)d(b)$. Then, $(A,\vdash,\dashv,\perp)$ turns to be a triangular algebra. Similarly, associative dialgebras (category denoted by \textsf{Dias.}), introduced
by J.-L. Loday in \cite{Loday}, are triangular algebras by setting $\perp=0$ since the two required operations $\vdash$ and $\dashv$ obey by definition the first 5 axioms.
\end{exam}
\begin{theo}
The operad \textit{Triang.} is dual in the sense of \cite{GK} to the operad \textit{3-Dend.}, that is 
\textit{Triang.}=\textit{3-Dend.}$^!$ and \textit{3-Dend.}=\textit{Triang.}$^!$.
\end{theo}
\Proof
We compute the dual of $3-Dend.$. Since our operads are regular, we know that the action of the symmetric group can be forgotten. We consider then only $\mathcal{P}_n$. The $K$-vector space generating
operations is $3-Dend._2:=K\prec \ \oplus \ K\succ \oplus \ K\bullet$. Set $OP:=\{\prec,\succ,\bullet\}$. The $K$-vector space
made out of three variables is $K[P \times P] \oplus K[P \times P]$. Its dimension is 18. The operad $3-Dend.$ is completely determined by some subspace $R \subset K[P \times P] \oplus K[P \times P]$. Denote by $(\circ_1)\circ_2$ (resp. $\circ_1(\circ_2))$, $\circ_i \in OP$, the basis vector of the first (resp. the second) summand $K[P \times P]$. Observe that $R$ is spanned by 6 vectors of the form $(\circ_1)\circ_2 - \circ_1(\circ_2)$ obtained from axioms of 3-dendriform algebras. Identify the dual of $K[P]$ with itself by identifying a basis vector with its dual. According to \cite{GK}, the dual operad $3-Dend.$ is then completely determined by $R^\perp \subset K[P \times P] \oplus K[P \times P]$, where $R^\perp$ is the orthogonal space of R under the quadratic form 
$
\begin{pmatrix}
 Id & 0\\
0 & -Id
\end{pmatrix}.
$
\noindent
Identify now $\prec$ to $\dashv$, $\succ$ to $\vdash$ and $\bullet$ to $\perp$. The $K$-vector space $R^\perp$ becomes the space $R^!$ spanned by the 12 vectors obtained from axioms of triangular algebras. For instance,
consider the vector $(\dashv)\dashv-\dashv(\dashv)$
of $R^!$, identified with $(\prec)\prec - \prec(\prec)$ and observe that for instance
$\bra (\prec)\prec - \prec(\prec); (\prec)\prec-\prec (\star) \ket = 1-1=0$ and so on.
\eproof
\Rk
The following functorial diagram,
\begin{center}
$
\begin{array}{ccccc}
\textbf{As.}& \longrightarrow &\textbf{Dias.}& \longrightarrow &\textbf{Triang.} \\
 & \searrow & \downarrow & \swarrow &  \\
& & \textbf{Leib.} & &
\end{array}
$

holds, where \textbf{Leib.} is the category of Leibniz algebras (a generalisation of Lie algebras \cite{Loday}).
\end{center}
\subsection{The free triangular algebra.}
Let $V$ be a $K$-vector space. 
Denote by $T(V)$ 
the tensor module, that is, $$T(V):=K \oplus V \oplus
V^{\otimes 2}\oplus \ldots \oplus V^{\otimes n} \oplus \ldots.$$
Through the paper, a tensor $v_1 \otimes \ldots \otimes v_p$ will be denoted sometimes, for commodity by $v_1, \ldots, v_p$ or by $v_1 \ldots v_p$, when no confusion is possible.
The algebraic object to consider is obviously $V\otimes T(V)^{\otimes 3}$. For esthetic reasons, we will work with an isomorphism copy written in an unusual way as:
\begin{center}
\begin{center}
$\bigtriangleup (V):=
\begin{array}{ccccc}
     &         & T(V) & & \\
    &         & \otimes &          &\\
T(V)& \otimes & V& \otimes & T(V)
\end{array}
$
\end{center} 
\end{center} 
Let $\psi:T(V)\longrightarrow K$ and $\Psi:T(V)^{\otimes 2}\longrightarrow K$ be the canonical projections. Define now three binary operations 
$\vdash,\dashv,\perp:\bigtriangleup^{\otimes 2} (V)\longrightarrow \bigtriangleup(V)$ as follows,
\begin{center}
$
\begin{array}{ccccccccccccccccccc}
 &  & H & &   & &    & & H' & & & & &  & \psi(H)H' & & \\
 &  & \otimes &  &  &\vdash & & & \otimes &  & &=&  &  & \otimes &  &\\
L & \otimes & v& \otimes & R & & L' & \otimes & v'& \otimes & R' & & LvRL' & \otimes & v'& \otimes & R',
\end{array}
$
\end{center} 
\begin{center}
$
\begin{array}{ccccccccccccccccccc}
 &  & H & &   & &    & & H' & & & & &  & H\psi(H') & & \\
 &  & \otimes &  &  &\dashv & & & \otimes &  & &=&  &  & \otimes &  &\\
L & \otimes & v& \otimes & R & & L' & \otimes & v'& \otimes & R' & & L & \otimes & v& \otimes & RL'vR',
\end{array}
$
\end{center} 
and
\begin{center}
$
\begin{array}{ccccccccccccccccccc}
 &  & H & &   & &    & & H' & & & & &  & \Psi(H,H')RL'v' & & \\
 &  & \otimes &  &  &\perp & & & \otimes &  & &=&  &  & \otimes &  &\\
L & \otimes & v& \otimes & R & & L' & \otimes & v'& \otimes & R' & & L & \otimes & v& \otimes & R',
\end{array}
$
\end{center}
for any $R \in V^{\otimes p_1}$, $L \in V^{\otimes p_2}$, $R \in V^{\otimes p_3}$, $R' \in V^{\otimes p_4}$, $L' \in V^{\otimes p_5}$, $H' \in V^{\otimes p_6}$ and any $v,v' \in V$ and extended by bilinearity then. 

\noindent
Let $V$ be a $K$-vector space. By definition, the free triangular algebra on $V$ is the triangular algebra $Triang.(V)$ equipped with a $K$-linear map
$i:V \hookrightarrow Triang.(V)$ such that for any $K$-linear map $f:V \longrightarrow T$, where $T$ is a triangular algebra over $K$, there exists a unique triangular algebra morphism $\phi$ turning the diagram,
\begin{center}
$
\begin{array}{ccc}
&i & \\
V& \hookrightarrow& Triang.(V) \\
  & & \\
& f \searrow  & \downarrow \phi \\
  & & \\
& & T 
\end{array}
$
\end{center}
commutative.
\begin{theo}
\label{freetria}
Let $V$ be a $K$-vector space.
The $K$-vector space $\bigtriangleup(V)$ equipped with the three operations just defined is the free triangular algebra on $V$.
\end{theo}
\Proof
Checking axioms of triangular algebras for $\bigtriangleup(V)$ is left to the reader. The map $i:V \hookrightarrow \bigtriangleup(V)$ is the composite:
\begin{center}
\begin{center}
$
V \simeq 
\begin{array}{ccccc}
     &         & K & & \\
    &         & \otimes &          &\\
K & \otimes & V& \otimes & K
\end{array}
\hookrightarrow \ \ \ \bigtriangleup (V).$
\end{center} 
\end{center} 
Let $f:V \longrightarrow T$ be a linear map, where $T$ is a triangular algebra. We construct $\phi$ as follows. First on monomials from $\bigtriangleup (V)$, second we extend it by $K$-linearity.
Therefore, define $\phi:\bigtriangleup (V) \longrightarrow T$ by:
%\begin{tiny}
\begin{center}
$\phi(X):= f(v_{-p}) \vdash \ldots \vdash f(v_{-1})\vdash [f(v_{0})\perp(f(w_{1}) \vdash \ldots \vdash f(w_{k}))]\dashv f(v_{1})\dashv \ldots \dashv f(v_{q}),$
\end{center}

\noindent
if,

$X:=\begin{array}{ccccc}
     &         & w_{1}\otimes \ldots \otimes w_{k} & & \\
    &         & \otimes &          &\\
 v_{-p}\otimes \ldots \otimes v_{-1} & \otimes  & v_0& \otimes &  v_{1}\otimes \ldots \otimes v_{q} 
\end{array}$

\noindent
and obviously by,

\begin{center}
$\phi(X):= f(v_{-p}) \vdash \ldots \vdash f(v_{-1})\vdash f(v_{0})\dashv f(v_{1})\dashv \ldots \dashv f(v_{q}),$
\end{center}

\noindent
if,

$X:=\begin{array}{ccccc}
     &         & 1 & & \\
    &         & \otimes &          &\\
 v_{-p}\otimes \ldots \otimes v_{-1} & \otimes  & v_0& \otimes &  v_{1}\otimes \ldots \otimes v_{q}. 
\end{array}$

\noindent
According to dimonoid calculus rules \cite{Loday}, the writings at the right hand sides do have a meaning. Observe that the bracket $[\ldots]$ can be dropped and is just used here to recall the symmetry shape of $\bigtriangleup (V)$. We prove now, that so defined,
$\phi$ is a morphism of triangular algebras. Let us start with the binary operation $\perp$. We replace, with a slight abuse of notation, tensors by capital letters so as to ease proofs.
On the one hand,

\begin{eqnarray*}
A&:=&\phi(\begin{array}{ccccc}
     &         & H & & \\
    &         & \otimes &  &\\
 L & \otimes  & v& \otimes &  R 
\end{array})\perp  \phi(\begin{array}{ccccc}
     &         & H' & & \\
    &         & \otimes &          &\\
 L' & \otimes  & v'& \otimes &  R' 
\end{array}):= \\
\\[0.5cm]
& =& (f(L)\vdash [f(v)\perp f(H)]\dashv f(R))\perp(f(L')\vdash [f(v')\perp f(H')]\dashv f(R')).
\end{eqnarray*}
Set $z:=f(L')\vdash [f(v')\perp f(H')]\dashv f(R')$, $y=[f(v)\perp f(H)]\dashv f(R)$, we get
$A= (f(L)\vdash y)\perp z = f(L)\vdash (y\perp z)$
\textit{via} Axiom 7. However, $y\perp z=([f(v)\perp f(H)]\dashv f(R))\perp z=(f(v)\perp [f(H)\dashv f(R)])\perp z =0$ \textit{via} Axiom 8 first and Axiom 9 then. The case $H=1$ and $H'\not=1$ give the same result and is left to the reader (apply Axioms 8, 7 and 10). If $H=H'=1$, then we get,
$A:=(f(L)\vdash f(v)\dashv f(R))\perp(f(L')\vdash f(v')\dashv f(R'))$. Set $z:=f(L')\vdash f(v')\dashv f(R')$, $y=f(v)\dashv f(R)$. Then, by applying Axiom 7 we get $(f(L)\vdash y)\perp z=
f(L)\vdash (y \perp z)$. However, $y \perp z:=(f(v)\dashv f(R)) \perp z= f(v)\perp (f(R) \vdash z)$ by Axiom 6. Moreover,
$f(R) \vdash z:=f(R) \vdash ((f(L')\vdash f(v'))\dashv f(R'))=(f(R) \vdash (f(L')\vdash f(v')))\dashv f(R')$ \textit{via} Axiom 3, which is equal to
$(f(R) \vdash f(L')\vdash f(v'))\dashv f(R')$ \textit{via} Axiom 5. Set $x:= f(R) \vdash f(L')\vdash f(v')$, we have $f(v) \perp (x \dashv f(R'))=(f(v) \perp x) \dashv f(R')$ \textit{via} Axiom 8. Summarising our computations, we find,
$$A:= f(L)\vdash(f(v)\perp (f(R) \vdash f(L')\vdash f(v')))\dashv f(R').$$ 
On the other hand, 
\begin{eqnarray*}
B&:=&\phi(\begin{array}{ccccc}
     &         & H & & \\
    &         & \otimes &  &\\
 L & \otimes  & v& \otimes &  R 
\end{array}\perp  \begin{array}{ccccc}
     &         & H' & & \\
    &         & \otimes &          &\\
 L' & \otimes  & v'& \otimes &  R' 
\end{array}):= 
\\[0.5cm]
& =& f(L)\vdash [f(v)\perp (f(R)\vdash(f(L')\vdash f(v'))] \dashv f(R').
\end{eqnarray*}
if $H=H'=1$ and vanishes otherwise.
Hence $A=B$.

\noindent
Concerning the binary operation $\vdash$, we get
on the one hand:
\begin{eqnarray*}
B&:=&\phi(\begin{array}{ccccc}
     &         & H & & \\
    &         & \otimes &  &\\
 L & \otimes  & v& \otimes &  R 
\end{array}\vdash  \begin{array}{ccccc}
     &         & H' & & \\
    &         & \otimes &          &\\
 L' & \otimes  & v'& \otimes &  R' 
\end{array}):= 
\\[0.5cm]
& =& f(L)\vdash f(v)\vdash f(R)\vdash f(L')\vdash [f(v') \perp f(H')]\dashv f(R'),
\end{eqnarray*}
if $H=1$ and vanishes otherwise. On the other hand,
\begin{eqnarray*}
A&:=&\phi(\begin{array}{ccccc}
     &         & H & & \\
    &         & \otimes &  &\\
 L & \otimes  & v& \otimes &  R 
\end{array}) \vdash  \phi(\begin{array}{ccccc}
     &         & H' & & \\
    &         & \otimes &          &\\
 L' & \otimes  & v'& \otimes &  R' 
\end{array}):= \\
\\[0.5cm]
& =& (f(L)\vdash [f(v)\perp f(H)]\dashv f(R))\vdash (f(L')\vdash [f(v')\perp f(H')]\dashv f(R')).
\end{eqnarray*}
Setting $z:=f(L')\vdash [f(v')\perp f(H')]\dashv f(R')$, we get by applying successively Axioms 3,4,5 and 11, 
\begin{eqnarray*}
A&=& ((f(L)\vdash [f(v)\perp f(H)])\dashv f(R))\vdash z = (f(L)\vdash [f(v)\perp f(H)]\vdash f(R))\vdash z \\
&=& (f(L)\vdash ([f(v)\perp f(H)]\vdash f(R)))\vdash z =0.
\end{eqnarray*}
If $H=1$, then it is easy to check the required equality. Therefore, $A=B$. Proceeding the same way
for the binary operation $\dashv$, shows that $\phi$ is a morphism of triangular algebras. Consequently, $\phi$ so constructed is unique since
it has to coincide with $f$ on $V$.
This completes the proof.
\eproof
\NB
As $Triang.$ is a regular binary and quadratic operad, the free triangular algebra over a $K$-vector space $V$ can be written as,
$$ Triang.(V):=\bigoplus_{n=1}^{\infty}Triang._n\otimes V^{\otimes n}.$$
On one generator, we get $Triang.(K):=\bigoplus_{n=1}^{\infty}Triang._n$. 
Set,
\begin{center}
\begin{center}
$\chi:=
\begin{array}{ccccc}
     &         & 1 & & \\
    &         & \otimes &          &\\
1& \otimes & v & \otimes & 1.
\end{array}$
\end{center} 
\end{center}
The $K$-vector space $(\bigtriangleup(K\chi),\vdash,\dashv,\perp)$ is the free triangular algebra over the generator $\chi$.
Observe that the dimension
of the degree $n$ part of $\bigtriangleup(K\chi):=\bigoplus_{n=1}^\infty Triang._n$
is the $n^{th}$ triangular number $\frac{n(n+1)}{2}$.
Let us describe in more details, the free triangular algebra on one generator. For all $n>0$, introduce the set $Tab(n):=\{
\begin{array}{ccc}
 &  p_3 & \\
\hline 
p_1& 1 & p_2 
\end{array}
; \ \ p_1+p_2+p_3=n-1
\}$.
Check that the cardinal of $Tab(n)$ is $\frac{n(n+1)}{2}$. We can turn the $K$-vector space $Tab_{\infty}:=\bigoplus_{n=1}^{\infty}KTab(n)$ into a triangular algebra by using the following bijection,
\begin{center}
$\begin{array}{ccccc}
     &         & v^{\otimes p_3} & & \\
    &         & \otimes &          &\\
v^{\otimes p_1}& \otimes & v & \otimes & v^{\otimes p_2}
\end{array} \longmapsto \ \ \ \begin{array}{ccc}
 &  p_3 & \\
\hline 
p_1& 1 & p_2 
\end{array},$
\end{center}
\noindent
where by convention $v^{\otimes 0}=1$. Consequently, $Tab_{\infty}:=\bigoplus_{n=1}^{\infty}KTab(n)$ is the free triangular algebra on the generator still denoted by $\chi =\begin{array}{ccc}
 &  0 & \\
\hline 
0& 1 & 0 
\end{array}$. 

\noindent
Now, let $\mathbb{P}^nHmg[X_0,X_1,X_2]$ be the $K$-vector space of homogeneous polynomials of degree $n-1$ over the commutative indeterminates $X_0$, $X_1$ and $X_2$. Observe that as a $K$-vector space, the dimension of $\mathbb{P}^n Hmg[X_0,X_1,X_2]$ is $\frac{n(n+1)}{2}$. For instance, $\mathbb{P}^3 Hmg[X_1,X_2,X_3]$ is spanned by $X_0^2$, $X_1^2$, $X_2^2$, $X_0X_1$, $X_0X_2$ and $X_1X_2$. Consider the bijection,
\begin{center}
$\begin{array}{ccccc}
     &         & v^{\otimes p_3} & & \\
    &         & \otimes &          &\\
v^{\otimes p_1}& \otimes & v & \otimes & v^{\otimes p_2}
\end{array} \longmapsto \ \ \ X_0^{p_1}X_1^{p_2}X_2^{p_3}.$
\end{center}
Set $\mathbb{P}^{\infty} Hmg[X_0,X_1,X_2]:= \bigoplus_{n=2}^{\infty} \mathbb{P}^n Hmg[X_0,X_1,X_2]$. This $K$-vector space inherits
 a triangular algebra structure \textit{via} the bijection
constructed above. Still denote by $\chi$ the image
of $\chi$ under this map. We get 
that,
$K\chi \oplus \mathbb{P}^{\infty} Hmg[X_0,X_1,X_2],$
is the free triangular algebra generated by $\chi$. Consequently, $\mathbb{P}^{\infty} Hmg[X_0,X_1,X_2]$ inherits  an operadic arithmetic, that is, a $K$-left-linear map (called usually the multiplication),
$$\circledast: \mathbb{P}^{\infty} Hmg[X_0,X_1,X_2] \otimes \mathbb{P}^{\infty} Hmg[X_0,X_1,X_2] \rightarrow \mathbb{P}^{\infty} Hmg[X_0,X_1,X_2],$$  distributive to the left as regards operations $\vdash, \ \dashv, \ \perp$ and consisting to replace the generator $\chi$ in the code describing
the left object by the right one. For instance,
$((\chi \perp \chi)\dashv \chi)\circledast z:=(z \perp z)\dashv z,$ where $z$ is a homogeneous polynomial of degree say $d$. The combinatorial object underlying the free triangular algebra
on one generator is no longer linear combinations of planar ternary trees but homogeneous polynomials, thus in relation to projective algebraic curves in the projective plane $\mathbb{P}^2(K)$. 

\noindent
We now relate these results to the construction of the homology of 3-dendriform algebras. Recall that
for all $n \in \mathbb{N}$, we defined $X_n:=\{(k,j); 0\leq k \leq j \leq n\}$.
The map $\eta: Tab(n) \longrightarrow X_n$ defined for all $n$ by,
\begin{center}
$
\begin{array}{ccc}
 &  p_3 & \\
\hline 
p_1& 1 & p_2 
\end{array} \longmapsto \ \ \begin{cases}(p_1,p_3) \ \textrm{if} \ p_1 \leq p_3, \\
(p_1,n-1-p_3) \ \textrm{if} \ p_1 > p_3,
\end{cases}$
\end{center}
is clearly a bijection and explain why we constructed the chain-complex of 3-dendriform algebras as we did, according to Ginzburg and Kapranov's results \cite{GK}.
\subsection{(Co)Homology of triangular algebras}
For any $t \in \treeg_n$, label its leaves from left to right starting from 0 to $2n$. Start with the leave 0, and begin to count from the place you reached. Every 3 leaves, place
the operation $\perp$ if the leave has the shape $\mid$, $\vdash$ if the shape is $/$ and $\dashv$ if the shape is $\setminus$. By convention the last leave remains unassigned. Proceeding that way, a tree from $\treeg_{n+1}$ will give $n$ binary operations, we label 1 to $n$ from left to right. For instance,
$\treeggd \in \treeg_3$ will give a $\dashv$ at leave 2 and $\vdash$ at leave 4. We have defined for $1\leq i\leq n$, a map,
$$ \circ_i: \treeg_{n+1} \longrightarrow \{\perp,\vdash,\dashv \},$$
assigning to each $i$ the corresponding operation
by the process just described. The image will be denoted by $\circ_i^t$. Denote by $\widetilde{\treeg}_n$ the set of trees of $\treeg_n$
whose leaves have been colored by the operations 
$\perp,\vdash,\dashv$ just as explained. For type-writing commodity, we will write the operations at the right hand side of the tree and the label of the leave it refers. We get a bijection
$tilde: \treeg_n \longrightarrow \widetilde{\treeg}_n$. For instance,
$$\treeggd \mapsto [\treeggd, (2;\dashv),(4;\vdash)], \ \ \treemmd \mapsto [\treemmd, (2;\perp),(4;\perp)].$$
For any $1\leq i \leq n$, define the face map,
$d_i:K\treeg_n \longrightarrow K\treeg_{n-1}$, on $\treeg_n$ first and extended by linearity 
then to be the composite $tilde^{-1} \ \circ del_i \ \circ \ tilde$, where the linear map,
$$del_i: K\widetilde{\treeg}_n \longrightarrow K\widetilde{\treeg}_{n-1},$$
assigns to a tree $t$, the tree $t'$ obtained from $t$ as follows. Localise the leave colored by the operation labelled by $i$. 
Remove the offspring of its father vertex if the three children are leaves, or if the middle leave is not colored. Otherwise, the result is zero.
For instance,
$d_2([\treeggd, (2;\dashv),(4;\vdash)]):=[\treegg, (2,\vdash)]$,  
$d_1([\treemmd, (2;\perp),(4;\perp)])=0$ since the middle leave (number 2) is colored and the offspring of its father vertex do not give only leaves, but
$d_2([\treemmd, (2;\perp),(4;\perp)])=[\treegm,(2;\perp)]$.
As well, $d_1([\treedmm, (2;\dashv),(4;\perp)])=0,$ since the father vertex associated with leave 2 has a middle children which is not a leave. However,
$d_2([\treedmm, (2;\dashv),(4;\perp)])=[\treedd,(2;\dashv)]$.
\noindent
Let $T$ be a triangular algebra over $K$. Define the module of $n$-chains by $C\treeg_n(T):=K[\treeg_n]\otimes T^{\otimes n}$.
Define a linear map $d: C\treeg_n(T) \longrightarrow C\treeg_{n-1}(T)$ by the following formula,
$$ d(t; v_1,\ldots, v_n):=\sum_{i=1}^{n-1} \ (-1)^{i+1}(d_i(t);v_1, \ldots, v_{i-1}, v_i \circ_i^t v_{i+1}, \ldots, v_n):=\sum_{i=1}^{n-1} \ (-1)^{i+1}d_i(t; v_1, \ldots, v_n),  $$
with a slight abuse of notation, and where $t\in \treeg_n$, $v_i \in T$.
\begin{prop}
The face maps $d_i:C\treeg_n(T) \longrightarrow C\treeg_{n-1}(T)$ satisfy the simplicial relations
$d_id_j=d_{j-1}d_i$ for any $1\leq i < j \leq n-1$.
Therefore $d\circ d=0$ and so $(C\treeg_*(T),d)$ is a chain-complex.
\end{prop}
\Proof
We check this idendity for the lowest dimension, that is,
$$d_1d_2=d_{1}d_1:C\treeg_3(T) \rightarrow C\treeg_{2}(T),$$ 
and 12 cases to detail.
\begin{itemize}
\item {Case $\treeddd$:}
$$ d_1d_2(\treeddd; x,y,z)=x\dashv (y \dashv z),$$
$$ d_1d_1(\treeddd; x,y,z)= (x \dashv y)\dashv z,$$
hence the equality by Axiom 1.
\item {Case $\treegdg$:}
$$ d_1d_2(\treegdg; x,y,z)=x\dashv (y \vdash z),$$
$$ d_1d_1(\treegdg; x,y,z)= (x \dashv y)\dashv z,$$
hence the equality by Axiom 2.
\item {Case $\treemgm$:}
$$ d_1d_2(\treemgm; x,y,z)=x\vdash (y \dashv z),$$
$$ d_1d_1(\treemgm; x,y,z)= (x \vdash y)\dashv z,$$
hence the equality by Axiom 3.
\item {Case $\treeggd$:}
$$ d_1d_2(\treeggd; x,y,z)=x\vdash (y \vdash z),$$
$$ d_1d_1(\treeggd; x,y,z)= (x \dashv y)\vdash z,$$
hence the equality by Axiom 4.
\item {Case $\treeggg$:}
$$ d_1d_2(\treeggg; x,y,z)=x\vdash (y \vdash z),$$
$$ d_1d_1(\treeggg; x,y,z)= (x \vdash y)\vdash z,$$
hence the equality by Axiom 5.
\item {Case $\treegmm$:}
$$ d_1d_2(\treegmm; x,y,z)=x\perp (y \vdash z),$$
$$ d_1d_1(\treegmm; x,y,z)= (x \dashv y)\perp z,$$
hence the equality by Axiom 6.
\item {Case $\treeGm$:}
$$ d_1d_2(\treeGm; x,y,z)=x\vdash (y \perp z),$$
$$ d_1d_1(\treeGm; x,y,z)= (x \vdash y)\perp z,$$
hence the equality by Axiom 7.
\item {Case $\treegmd$:}
$$ d_1d_2(\treegmd; x,y,z)=x\perp (y \dashv z),$$
$$ d_1d_1(\treegmd; x,y,z)= (x \perp y) \dashv z,$$
hence the equality by Axiom 8.
\item {Case $\treemmg$:}
$$ d_1d_2(\treemmg; x,y,z)=0,$$
$$ d_1d_1(\treemmg; x,y,z)= (x \perp y)\perp z,$$
hence the equality by Axiom 9.
\item {Case $\treemmd$:}
$$ d_1d_2(\treemmd; x,y,z)=x\perp (y \perp z),$$
$$ d_1d_1(\treemmd; x,y,z)= 0,$$
hence the equality by Axiom 10.
\item {Case $\treeggm$:}
$$ d_1d_2(\treeggm; x,y,z)=0,$$
$$ d_1d_1(\treeggm; x,y,z)= (x \perp y)\vdash z,$$
hence the equality by Axiom 11.
\item {Case $\treedmm$:}
$$ d_1d_2(\treedmm; x,y,z)=x\dashv (y \perp z),$$
$$ d_1d_1(\treedmm; x,y,z)= 0,$$
hence the equality by Axiom 12.
\end{itemize} 
The general case splits into two cases. If $j=i+1$, then the proof follows from the low dimension cases and from axioms of triangular algebras. The case
$j>i+1$ is straightforward.
\eproof

\noindent
We get a chain-complex,
$$C\treeg_*(T):\ \ \ldots \rightarrow K[\treeg_n]\otimes T^{\otimes n} \rightarrow \ldots 
\rightarrow  K[\treeg_3]\otimes T^{\otimes n} \rightarrow K[\treeg_2]\otimes T^{\otimes n} \xrightarrow{\perp,\vdash,\dashv} T.$$
This allows us to define the homology of a triangular algebra $T$ as the homology of the chain-complex $C\treeg_*(T)$, that is $H\treeg_n(T):=H_n(C\treeg_*(T),d)$, $n>0$. The cohomology of $T$ is by definition $H\treeg^n(T):=H^n(C\treeg_*(T),K)$, $n>0$.
For $\bigtriangleup(V)$, the free triangular algebra over the $K$-vector space $V$, we get 
$H\treeg_1(\bigtriangleup(V))\simeq V$, since
by definition, 
$$H\treeg_1(\bigtriangleup(V)):= \bigtriangleup(V)/ \{x\perp y, x\vdash y, x\dashv y; x,y \in \bigtriangleup(V) \} \simeq \begin{array}{ccccc}
     &         & K & & \\
    &         & \otimes &          &\\
K& \otimes & V & \otimes & K
\end{array} \simeq V.$$
We conjecture that for $n>1$,     $H\treeg_n(\bigtriangleup(V))=0$, that is the operad $Triang.$ is Koszul.
\section{On planar rooted $m$-ary trees and tetrahedral numbers}
Fix an integer $m>1$.
In this section, we will generalise our results to 
the planar rooted $m$-ary trees, $m$-trees for short and their relations to tetrahedral numbers.
\subsection{The operad $m-Dend.$}
\begin{defi}{}
A $K$-vector space $T$ is a $m$-dendriform algebra if it is equipped with $m$ binary operations
$\prec,\succ,\bullet_2, \ldots, \bullet_{m-1}: T^{\otimes 2} \longrightarrow T$ verifying for all $x,y,z \in T,$ and for all $2 \leq i \leq m-1$,
the $\frac{m(m+1)}{2}$ axioms.
$$ (x\prec y)\prec z = x \prec (y \star z), \ \  \  \  \ (x \prec y)\bullet_i z =x\bullet_i (y \succ z) $$ 
$$ (x \succ y)\prec z = x \succ (y \prec z), \ \  \  \  \ (x \succ y)\bullet_i z =x\succ (y \bullet_i z) $$
$$(x \star y)\succ z =x\succ (y \succ z) , \ \  \  \  \ (x\bullet_i  y)\prec z = x \bullet_i (y \prec z), $$
where $x \star y := x\prec y + x\succ y$ and,
$$(x\bullet_i  y)\bullet_j z = x \bullet_i (y \bullet_j z), $$
for all $2 \leq i<j \leq m-1$. 
A $m$-dendriform algebra is said to be involutive if it is equipped with an involution $\dagger$
verifying equalities written in (\ref{inv}).
\end{defi}
\Rk \textbf{(Opposite $m$-dendriform algebra.)}
Let $T$ be a $m$-dendriform algebra. Define new operations by:
$$x \prec' y:= y \succ x; \ x \succ' y:= y \prec x; \ x \bullet_i' y:= y \bullet_{m+1-i} x,$$
for all $2 \leq i \leq m-1$. Then, the $K$-vector space $T$ equipped with these operations is a new
$m$-dendriform algebra, denoted by $T^{op}$, called the opposite $m$-dendriform algebra. A $m$-dendriform algebra is said to be commutative if
$T^{op}=T$.

\noindent
Therefore, each $m$ gives birth to a regular, binary and quadratic operad denoted by $m-Dend.$, and a category denoted by $m$-\textsf{Dend.}.
The functorial diagram between categories holds:
\begin{center}
$
\begin{array}{cccc} 
2-\textsf{Dend.} \longrightarrow &     3-\textsf{Dend.} \longrightarrow \ldots &\longrightarrow m-\textsf{Dend.} & \ldots \longrightarrow \\ 
 _+ \searrow &  \downarrow \ _+ &  \swarrow _+ & \\
&    \textsf{As.} &  & \\ 
\end{array} 
$
\end{center}
We recall some general facts about $m$-ary trees.
Their names comes from the fact that each father vertex has exactly $m$ children. The degree of such a tree will be the number of internal vertices. By 
$\overset{_m}{\treeg}_n$, we mean the set of $m$-ary trees of degree $n$. It is known that the cardinal of 
$\overset{_m}{\treeg}_n$ is $\frac{(mn)!}{n!(n(m-1)+1)!}$.
Each $m$-tree can be decomposed in a unique way \textit{via} the so called grafting operation defined for all $p_1, \ldots, p_m \in \mathbb{N}$ as follows, 
$$\vee:\overset{_m}{\treeg}_{p_1} \times \ldots \times \overset{_m}{\treeg}_{p_m}  \longrightarrow \overset{_m}{\treeg}_{_{\sum_{i=1}^m p_i+1}},$$
$$ (t_1,\ldots,t_m) \longmapsto t_1\vee \ldots \vee t_m,$$
the last symbol meaning that the roots of $t_i$ are glued together and that a new root is created. The tree $\overset{_m}{\treeg}$ represents by definition
the tree $\mid \vee \ldots \vee  \mid$ $m$-times and is called a $m$-corolla. For instance, $\overset{_3}{\treeg}:=\treeg$, $\overset{_4}{\treeg}:=\treeCor$, and so forth.
There exists also an involution still denoted by $\dagger:\overset{_m}{\treeg}_n \mapsto \overset{_m}{\treeg}_n$ for all $n$, defined recursively by $t^\dagger=t_m^\dagger\vee \ldots \vee t_1^\dagger$ if $t=t_1\vee \ldots \vee t_m$. Pictorially, this is the mirror image \textit{via} the central axis of the tree. 
Over the $K$-vector space $K[\overset{_m}{\treeg}_\infty]:=K\mid \oplus \widehat{K[\overset{_m}{\treeg}_\infty]}$, where $\widehat{K[\overset{_m}{\treeg}_\infty]}:=\bigoplus_{n>0} K\overset{_m}{\treeg}_n$, we introduce recursively the following binary operations first on the trees. They are naturally extended by bilinearity to the whole $\widehat{K[\overset{_m}{\treeg}_\infty]}$. Let $p,q>0$ and
set for any $t=t_1\vee \ldots \vee t_m \in \overset{_m}{\treeg}_p$ and $r=r_1\vee \ldots \vee r_m \in \overset{_m}{\treeg}_q$,
\begin{eqnarray}
t\succ r &=& (t\star r_1)\vee \ldots \vee r_m, \\ 
t\prec r &=& t_1\vee \ldots \vee (t_m\star r),\\
t\bullet_2 r &=& r_1\vee (r_2 \vee \ldots \vee r_m \star t_1 \vee t_2) \vee \ldots \vee t_m,\\
t\bullet_3 r &=& r_1\vee r_2\vee (r_3 \vee \ldots \vee r_m \star t_1 \vee t_2\vee t_3) \vee \ldots \vee t_m,\\
\vdots \ \ \ &=& \ \ \ \vdots \hskip3cm \vdots\\
t\bullet_{_{m-1}} r &=& r_1\vee \ldots \vee (r_{m-1} \vee r_m \star t_1 \vee \ldots \vee t_{m-1}) \vee  t_m,\\
\mid \prec t &=& 0= t\succ \mid,\hskip1cm
t \prec \mid = t= \mid \succ t,\\
\mid\bullet_i t &=& 0=t\bullet_i \mid,
\end{eqnarray} 
where as usual $x \star y= x\prec y + x \succ y$ is
by construction associative. The symbols $1\prec 1$, $1\succ 1$ and $1 \bullet_i 1$ are not defined. Observe that $
\mid \star t = t= t\star \mid$ and that our three operations respect the natural grading of $\widehat{K[\overset{_m}{\treeg}_\infty]}$ since $\bullet_i,\succ,\prec: K\overset{_m}{\treeg}_p \otimes K\overset{_m}{\treeg}_q \longrightarrow K\overset{_m}{\treeg}_{p+q}$.
As expected, we get the following results whose proofs are similar of those decribed for $m=3$.
\begin{theo}
Fix an integer $ m \geq 2$.
\begin{itemize}
\item {Equipped with these $m$ binary operations, $\widehat{K[\treeg_\infty]}$ is the free $m$-dendriform algebra generated by $\overset{_m}{\treeg}$ which is involutive when equipped with the involution $\dagger$.}
\item {Extend the involution $\dagger$ as explained in the first section. Then, there exists a structure of involutive and connected $\mathcal{P}$-Hopf algebra on  $K[\overset{_m}{\treeg}_\infty]:=K\mid \oplus \widehat{K[\overset{_m}{\treeg}_\infty]}$.}
\item {Let $V$ be a $K$-vector space.
The following operations turn
$\bigoplus_{n>0} K[\overset{_m}{\treeg}_n]\otimes V^{\otimes n}$ into a $m$-dendriform algebra:
\begin{eqnarray*}
t\otimes \omega \prec t'\otimes \omega' &:=& t\prec t'\otimes \omega\omega', \\
t\otimes \omega \succ t'\otimes \omega' &:=& t\succ t'\otimes \omega\omega', \\
t\otimes \omega \bullet_i t'\otimes \omega' &:=& t\bullet_i t'\otimes \omega\omega'. \\
\end{eqnarray*}
Therefore, the unique $m$-dendriform algebra map $m-Dend.(V) \longrightarrow \bigoplus_{n>0} K[\overset{_m}{\treeg}_n]\otimes V^{\otimes n},$ which sends the generator $v\in V$ to $\overset{_m}{\treeg} \otimes v$ is an isomorphism.}
\end{itemize} 
\end{theo}
\subsection{A possible link to genomics}
In general, mathematicians model genomics
of species over the alphabet $\mathbb{A}=\{A,C,G,T\}$, a genomic sequence being an element of the dictionary $\mathbb{A}^*:=\cup_{n=0}^{\infty} \mathbb{A}^n$, where $\mathbb{A}^0$ contains the empty sequence and $\mathbb{A}^n$ contains all sequences of length $n$. Another possibility is to consider rooted planar 4-ary trees as in \cite{Dimitriadn}. We go further on that idea. Label the four leaves of $\overset{_4}{\treeg}:=\treeCor$ from left to right by $A$, $C$, $G$ and $T$. Concatening two letters
$A$ and $G$ for instance, giving $AC$, means
in the language of trees to consider a stacking of three trees as follows.
\begin{center}
\includegraphics*[width=2cm]{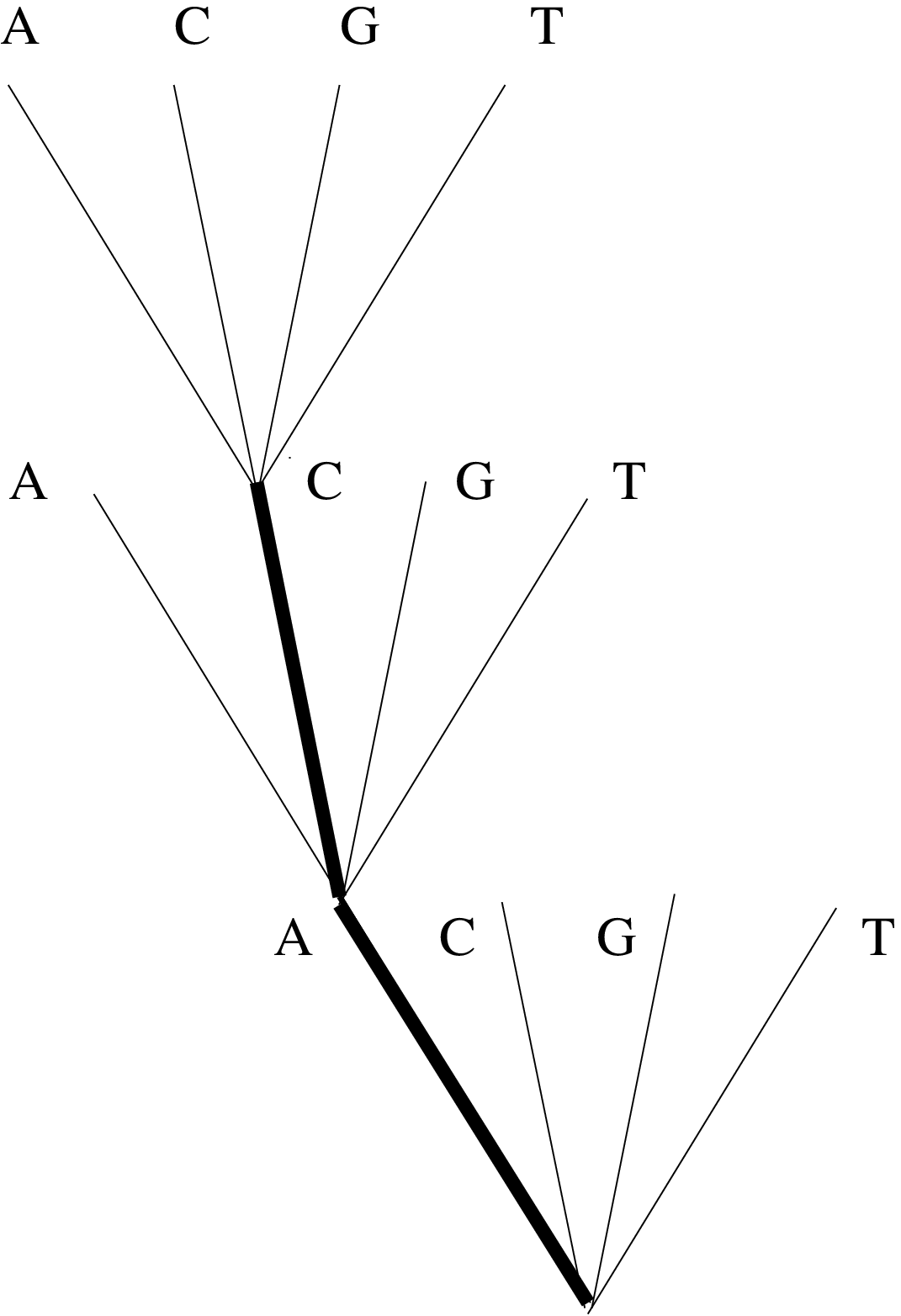}
\end{center}
We go first to $A$, then to $C$ and after the possibilities are open. In genomics, we have a pairing $A \leftrightarrow T$ and $C \leftrightarrow G$ which appears in the structure of DNA chains. In the language of trees, we have an involution $\dagger$ which takes into account exactly that pairing. Mutations can also be taken into account since one can use other branches of the tree $\treeCor$ to concatenate letters. If for instance, instead of considering $AC$ in our example,
we consider $AC$ and $AT$ (the letter $C$ mutates in $T$. We can consider the tree above and adding on it another tree to obtain:
\begin{center}
\includegraphics*[width=2cm]{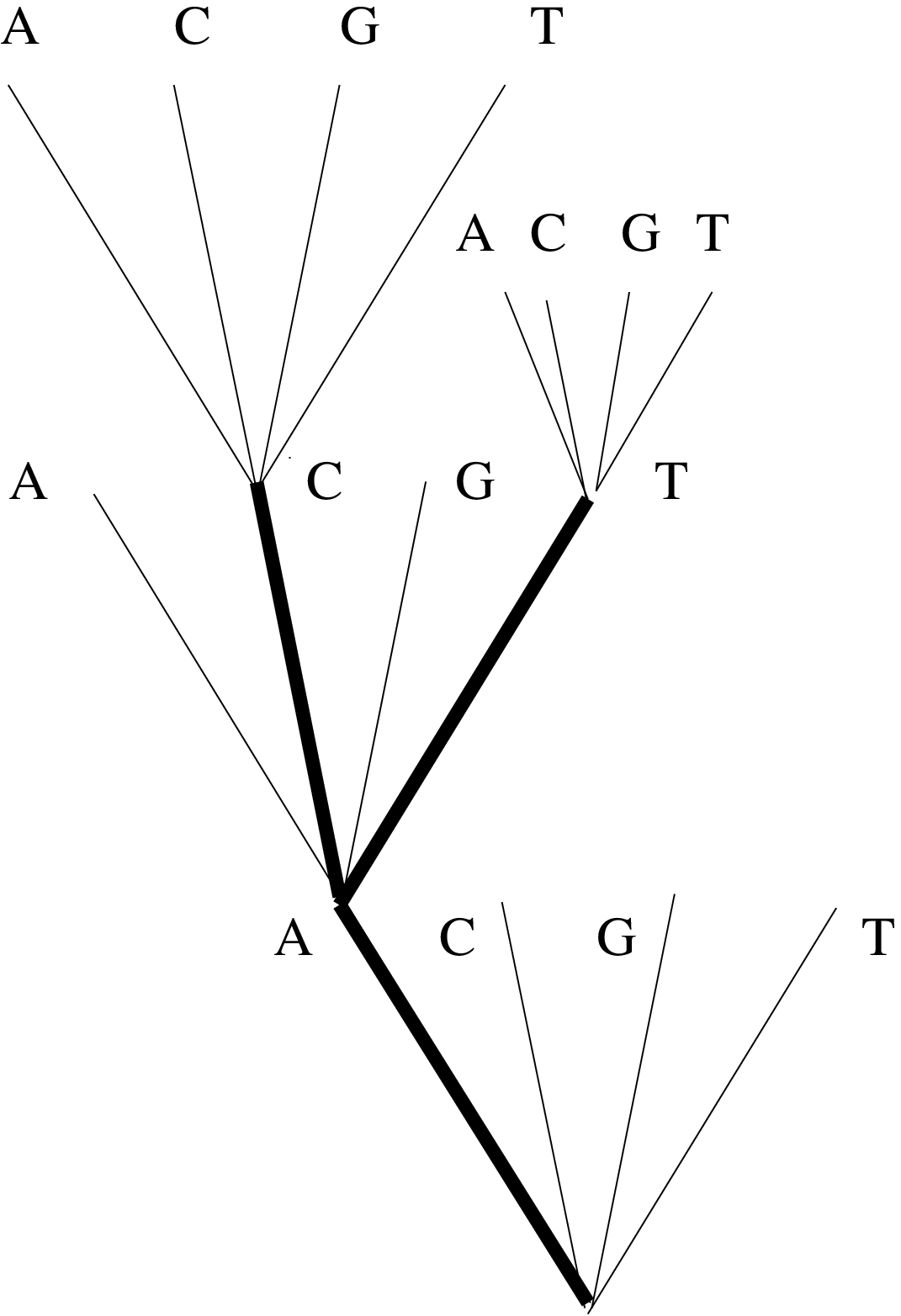}
\end{center}
Consequently, one sees that instead of considering
the set $\mathbb{A}^*$, one could consider the free
4-dendriform algebra over the generator $\treeCor$ equipped with the involution $\dagger$. This algebraic structure, as we saw, is much richer that the simple set $\mathbb{A}^*$ and instead of considering four symbols, just one is enough.
\subsection{The operad $m-Tetra.$}
We have seen that the operad $3-Dend.$ and $Triang.$ were dual in the sense of Ginzburg and Kapranov. The r\^ole of triangular numbers
will be played by the tetrahedral numbers of dimension $k$, where $k=m-1$. The $n^{th}$ tetrahedral number
of dimension $k$ is by definition the number
$t_{_{[k]}}^{[n]}:=\frac{n(n+1) \ldots (n+k-1)}{k!}$ (see N.J.A. Sloane Online Encyclopedia of Integers for instance).
For instance, for $k=3$, we got the triangular numbers $1,3,6,10,15, \ldots$ present in the definition of the Poincar\'e series of $f_{Triang.}$ in relation to $1,3,12,56, \ldots$ present in the definition of the Poincar\'e series of $f_{3-Dend.}:=f_{Triang.^!}$. Observe, for instance by using Maple, that these two series are inverse each other for the usual composition of functions. Instead of dimension 2 (triangles are drawn in a plane), consider dimension 3. Triangles become tetraedrons, hence the tetrahedral numbers.
They are $1,4,10,20,\ldots$ and they will be present in the definition of the Poincar\'e series of the operad $Tetra.$ defined below. This sequence  represents the coefficients of the inverse function  (for the usual composition) of the Poincar\'e series whose coefficients are $1,4,22,\ldots$ counting the numbers of 4-ary trees of degree $n$. In addition, concerning Poincar\'e series, we obtain that $f_{4-Dend.}$ and $f_{Tetra.}$ are inverse each other for the usual composition of functions. If we add another dimension, we will get the tetrahedral numbers of dimension 4: $1,5,15, \ldots$ related as expected to
$1,5,35, \ldots$ counting the numbers of 5-ary trees of degree $n$ with the same conclusions and so forth.
We now define these types of algebras.
\begin{defi}{}
\label{deftetra}
A $m$-tetrahedral algebra $T$ is a $K$-vector space equipped with $m$ binary operations $\vdash,\dashv,\perp_2, \ldots,\perp_{_{m-1}}: T^{\otimes 2} \longrightarrow
T$ verifying for all $x,y,z \in T$, for all $2 \leq i\leq m-1$, the following $\frac{n(3n-1)}{2}$ axioms.

\begin{eqnarray}
\begin{cases}
  (x\dashv y)\dashv z= x\dashv(y\dashv z),& 
  (x\dashv y)\perp_i z=x \perp_i(y\vdash z),\\
  (x\dashv y)\dashv z= x\dashv(y\vdash z),& 
  (x\vdash y)\perp_i z=x \vdash(y\perp_i z),\\
  (x\vdash y)\dashv z= x\vdash(y\dashv z),& 
  (x\perp_i y)\dashv z=x \perp_i (y\dashv z),\\
  (x\dashv y)\vdash z= x\vdash(y\vdash z),& 
  (x\perp_i y)\perp_i z=0 =  x \perp_i(y\perp_i z),\\
  (x\vdash y)\vdash z= x\vdash(y\vdash z),& 
  (x\perp_i y)\vdash z=0 = x \dashv(y\perp_i z),
\end{cases}
\end{eqnarray}
and also for all $2 \leq i<j\leq m-1$,
\begin{eqnarray}
\label{LL}
\begin{cases}
(x\perp_i y)\perp_j z = x \perp_i(y\perp_j z),\\
(x\perp_j y)\perp_i z=0 =  x \perp_j(y\perp_i z).
\end{cases}
\end{eqnarray}
\end{defi}
Any associative algebra $(A, \cdot)$ is a $m$-tetrahedral algebra by setting $\vdash=\cdot=\dashv$ and $\perp_i=0$ for all $2\leq i \leq m-1$.
Axioms above give birth to regular, binary and quadratic operads denoted by $m-Tetra.$ and categories denoted by $m$-\textsf{Tetra.}.
For convenience, we set $3-Tetra.:=Triang.$, $3$-\textsf{Tetra.}:=\textsf{Triang.} and $4-Tetra.:=Tetra.$ and $4$-\textsf{Tetra.}:=\textsf{Tetra.}.
\noindent
The following functorial diagram holds:
\begin{center}
$
\begin{array}{cccc}
\textbf{As.}\rightarrow  &\textbf{Dias.}\rightarrow &\textbf{Triang.} \rightarrow & \textbf{Tetra.} \rightarrow \ldots \\
  \searrow & \downarrow & \swarrow & \swarrow \\
& \textbf{Leib.}&    &
\end{array}
$
\end{center}
\Rk \textbf{(Opposite $m$-tetrahedral algebra.)}
Let $T$ be a $m$-tetrahedral algebra. Define new operations by:
$$x \vdash' y:= y \dashv x; \ x \dashv' y:= y \vdash x; \ x \perp_i' y:= y \perp_{m+1-i} x,$$
for all $2 \leq i \leq m-1$. Then, the $K$-vector space $T$ equipped with these operations is a new
$m$-tetrahedral algebra, denoted by $T^{op}$, called the opposite $m$-tetrahedral algebra. A $m$-tetrahedral algebra is said to be commutative if
 $T^{op}=T$.
\Rk
For any $m$-tetrahedral algebra $T$, let $as(T)$
be the quotient of $T$ by the ideal generated by the elements $x \vdash y - x\dashv y$ and $x \perp_i y$, for all $2 \leq i \leq m-1$ and $x,y \in T$. Observe that $as(T)$ is an associative algebra and that the functor $as.(-):$\textsf{$m$-Tetra.}$\rightarrow$ \textsf{As.} is left adjoint to the functor $inc:$ \textsf{As.} $\rightarrow$ \textsf{$m$-Tetra.}.

\noindent
As expected, we get the following result. 
\begin{theo}
For each $m \geq 2$,
the operad $m$-\textit{Tetra.} is dual in the sense of \cite{GK} to the operad \textit{$m$-Dend.}, that is 
$m$-\textit{Tetra.}=\textit{$m$-Dend.}$^!$ and \textit{$m$-Dend.}=$m$-\textit{Tetra.}$^!$.
\end{theo}
\Proof
Straightforward.
\eproof

\noindent
Let $V$ be a $K$-vector space and $T(V):= K \oplus V \oplus V^{\otimes 2} \oplus \ldots \oplus V^{\otimes n}  \oplus \ldots$. Let the linear map, $$\Psi: \underset{m-1 \ times}{\underbrace{T(V)\otimes \ldots \otimes T(V)}} \longrightarrow K,$$ 
and the linear map,
$$\psi: \underset{m-2 \ times}{\underbrace{T(V)\otimes \ldots \otimes T(V)}} \longrightarrow K,$$
be the canonical projections.
By definition, the free $m$-tetrahedral algebra on $V$ is the $m$-tetrahedral algebra $m-Tetra.(V)$ equipped with a $K$-linear map
$i:V \hookrightarrow m-Tetra.(V)$ such that for any $K$-linear map $f:V \longrightarrow T$, where $T$ is a $m$-tetrahedral algebra over $K$, there exists a unique $m$-tetrahedral algebra morphism $\phi$ turning the diagram,
\begin{center}
$
\begin{array}{ccc}
&i & \\
V& \hookrightarrow& m-Tetra.(V) \\
  & & \\
& f \searrow  & \downarrow \phi \\
  & & \\
& & T 
\end{array}
$
\end{center}
commutative.
\begin{theo}
Let $V$ be a $K$-vector space. Consider the $K$-vector space, 
$$m-Tetra(V):=T(V)\otimes [V\otimes T(V) \otimes \ldots \otimes T(V)] \otimes T(V).$$ 
Then, equipped with operations, $\perp_i$ $2 \leq i \leq m-1$, defined by,
\begin{small}
\begin{eqnarray*}
\omega_1\otimes [v\otimes \omega_2 \otimes \ldots \otimes \omega_{_{m-1}}]&\otimes &\omega_{m} \perp_2  \omega'_1\otimes [v'\otimes \omega'_2 \otimes \ldots \otimes \omega'_{_{m-1}}]\otimes \omega'_{m} \\ &=& \Psi(\omega_2,\ldots ,\omega_{m-1},\omega'_2)\omega_1\otimes [v\otimes \omega_m\omega'_1v'\otimes \omega'_3 \otimes \ldots \otimes \omega'_{m-1}]\otimes \omega'_{m},\\
\omega_1\otimes [v\otimes \omega_2 \otimes \ldots \otimes \omega_{_{m-1}}]&\otimes &\omega_{m} \perp_3  \omega'_1\otimes [v'\otimes \omega'_2 \otimes \ldots \otimes \omega'_{_{m-1}}]\otimes \omega'_{m} \\ &=& \Psi(\omega_3,\ldots ,\omega_{m-1},\omega'_2,\omega'_3)\omega_1\otimes [v\otimes \omega_2 \otimes \omega_m\omega'_1v'\otimes \omega'_4 \otimes \ldots \otimes \omega'_{m-1}]\otimes \omega'_{m},\\
\omega_1\otimes [v\otimes \omega_2 \otimes \ldots \otimes \omega_{_{m-1}}]&\otimes &\omega_{m}  \perp_4  \omega'_1\otimes [v'\otimes \omega'_2 \otimes \ldots \otimes \omega'_{_{m-1}}]\otimes \omega'_{m} \\ &=& \Psi(\omega_4,\ldots ,\omega_{m-1},\omega'_2,\omega'_3,\omega'_4)\omega_1\otimes [v\otimes \omega_2 \otimes \omega_3 \otimes \omega_m\omega'_1v'\otimes \omega'_5 \otimes \ldots \otimes \omega'_{m-1}]\otimes \omega'_{m},\\
\vdots \ \ \ &=& \ \ \ \vdots \hskip3cm \vdots\\
\omega_1\otimes [v\otimes \omega_2 \otimes \ldots \otimes \omega_{_{m-1}}]&\otimes &\omega_{m}  \perp_{_{m-1}}   \omega'_1\otimes [v'\otimes \omega'_2 \otimes \ldots \otimes \omega'_{_{m-1}}]\otimes \omega'_{m} \\ &=& \Psi(\omega_{m-1},\omega'_2,\ldots,\omega'_{m-1})\omega_1\otimes [v\otimes \omega_2 \otimes \ldots \otimes \omega_{m-2}\otimes \omega_m\omega'_1v']\otimes \omega'_{m},\\
\end{eqnarray*}
and,
\begin{eqnarray*}
\omega_1\otimes [v\otimes \omega_2 \otimes \ldots \otimes \omega_{_{m-1}}]\otimes \omega_{m} &\dashv &  \omega'_1 \otimes [v'\otimes \omega'_2 \otimes \ldots \otimes \omega'_{_{m-1}}]\otimes \omega'_{m}  \\ &=& \psi(\omega'_2, \ldots, \omega'_{_{m-1}})\omega_1\otimes [v\otimes \omega_2 \otimes \ldots \otimes \omega_{_{m-1}}]\otimes \omega_{m}\omega'_1v'\omega'_{m},\\
\omega_1\otimes [v\otimes \omega_2 \otimes \ldots \otimes \omega_{_{m-1}}]\otimes \omega_{m} &\vdash &  \omega'_1 \otimes [v'\otimes \omega'_2 \otimes \ldots \otimes \omega'_{_{m-1}}]\otimes \omega'_{m}  \\ &=& \psi(\omega_2, \ldots, \omega_{_{m-1}})\omega_1 v \omega_{m}\omega'_1 \otimes [v'\otimes \omega'_2 \otimes \ldots \otimes \omega'_{_{m-1}}]\otimes \omega'_{m}.
\end{eqnarray*}
\end{small}
$m-Tetra.(V)$ is the free $m$-tetrahedral algebra
on $V$. Consequently, the Poincar\'e series associated with the operad $m-Tetra$ is,
$$ f_{m-tetra}(x):=\sum_{n\geq 1} (-1)^n t_{_{[m-1]}}^{[n]}x^n = \frac{-x}{(1+x)^m}.$$
\end{theo}
\NB
Observe that the brackets $[\ldots]$ can be dropped
and are here only to recall the symmetry in the definitions of operations $\vdash$ and $\dashv$.
\Proof
The proof is similar to the one of Theorem \ref{freetria} and will be only sketched to avoid tedious computations. 
Therefore, checking axioms of $m$-tetrahedral algebras is left to the reader. The map $i:V \hookrightarrow m-Tetra.(V)$ is the composite:
$$ V \simeq K\otimes [V \otimes K\otimes \ldots \otimes K] \otimes K \hookrightarrow \ \ m-Tetra(V). $$
Let $f:V \longrightarrow T$ be a linear map, where $T$ is a $m$-tetrahedral algebra. We construct $\phi$ as follows. First on monomials from $m-Tetra(V)$, then extended by $K$-linearity.
Therefore, define $\phi:m-Tetra(V) \longrightarrow T$ by:
\begin{footnotesize}
\begin{eqnarray*}
\phi(X):= 
f(v_{-p}) \vdash \ldots \vdash f(v_{-1})\vdash [f(v_{0})&\perp_2 & (f(^2w_{1}) \vdash \ldots \vdash f(^2w_{k_2}))\\
&\perp_3 &(f(^3w_{1}) \vdash \ldots \vdash f(^3w_{k_3}))\\
 & \vdots & \\
&\perp_{_{m-1}} & (f(^{(m-1)}w_{1}) \vdash \ldots \vdash f(^{(m-1)}w_{k_{m-1}}))]\dashv f(v_{1})\dashv \ldots \dashv f(v_{q}),
\end{eqnarray*}
\end{footnotesize}
if,
\begin{footnotesize}
\begin{eqnarray*}
X:= 
v_{-p} \vdash \ldots \vdash v_{-1}\vdash [v_{0}&\perp_2 & (^2w_{1} \vdash \ldots \vdash ^2w_{k_2})\\
&\perp_3 &(^3w_{1} \vdash \ldots \vdash ^3w_{k_3})\\
 & \vdots & \\
&\perp_{_{m-1}} & (^{(m-1)}w_{1} \vdash \ldots \vdash ^{(m-1)}w_{k_{m-1}})]\dashv v_{1}\dashv \ldots \dashv v_{q},
\end{eqnarray*}
\end{footnotesize}
where the $v_j$ and the $^iw_j$ belong to $V$ and where it is understood that the non-existence of a line $\perp_i \ldots$, $2 \leq i \leq m-1$, in the definition of $X$ entails
the vanishing of the line $\perp_i \ldots$ in the 
construction of $\phi$.
According to dimonoid calculus rules \cite{Loday} and the first group of equations (\ref{LL}) of Definition \ref{deftetra}, the writings at the right hand side of the equalities do have a meaning. The proof that
$\phi$ is a morphism of $m$-tetrahedral algebras
is tedious but not difficult and is left to the reader. The unicity of such a morphism is now straightforward since it has to coincide on $V$ with $f$.
The last claim concerning the Poincar\'e series associated with the operad $m-Tetra$ is easy to compute.
\eproof

\subsection{Relations with homogeneous polynomials}
Set $\chi:= 1 \otimes [v \otimes 1 \otimes \ldots \otimes 1 ] \otimes 1$. 
Let $\mathbb{P}^nHmg[X_0,X_1,\ldots, X_{m-1}]$ be the $K$-vector space of homogeneous polynomials of degree $n-1$ over the commutative indeterminates $X_0,X_1,\ldots, X_{m-1}$. Observe that as a $K$-vector, the dimension of $\mathbb{P}^n Hmg[X_0,X_1,\ldots, X_{m-1}]$ is $t_{_{[m-1]}}^{[n]}$. Consider the bijection,
$$v^{\otimes p_1} \otimes [v \otimes v^{\otimes p_2} \otimes \ldots \otimes v^{\otimes p_{m}} ] \otimes 1 \longmapsto \ \ \ X_0^{p_1}X_1^{p_2}\ldots X_{m-1}^{p_m}.$$
Set $\mathbb{P}^{\infty} Hmg[X_0,X_1,\ldots, X_{m-1}]:= \bigoplus_{n=2}^{\infty} \mathbb{P}^n Hmg[X_0,X_1,\ldots, X_{m-1}]$. This $K$-vector space inherits
 a $m$-tetrahedral algebra structure \textit{via} the bijection
constructed above. Still denote by $\chi$ the image
of $\chi$ under this map. We get 
that,
$K\chi \oplus \mathbb{P}^{\infty} Hmg[X_0,X_1,\ldots, X_{m-1}],$
is the free $m$-tetrahedral algebra generated by $\chi$. Consequently, $\mathbb{P}^{\infty} Hmg[X_0,X_1,\ldots, X_{m-1}]$ inherits an operadic arithmetic, that is, a $K$-left-linear map (called usually the multiplication),
$$\circledast: \mathbb{P}^{\infty} Hmg[X_0,X_1,\ldots, X_{m-1}]^{\otimes 2} \rightarrow \mathbb{P}^{\infty} Hmg[X_0,X_1,\ldots, X_{m-1}],$$  distributive to the left as regards operations $\vdash, \ \dashv, \ \perp_i$ and consisting to replace the generator $\chi$ in the code describing
the left object by the right one. For instance,
$((\chi \perp_i \chi)\dashv \chi)\circledast z:=(z \perp_i z)\dashv z,$ where $z$ is a homogeneous polynomial of degree say $d$. The combinatorial object underlying the free $m$-tetrahedral algebra
on one generator is no longer linear combinations of planar $m$-ary trees but homogeneous polynomials over $m$ commutative indeterminates, thus related to projective algebraic hypersurfaces in the projective space $\mathbb{P}^{m-1}(K)$. 
\subsection{(Co)Homology of $m$-tetrahedral algebras}
For any $t \in \overset{_m}{\treeg}_n$, label its leaves from left to right starting from 0 to $(m-1)n$. Start with the leave 0, and begin to count from the place you reached. Every $m$ leaves, place
the operation $\perp_{m+1-i}$, $2\leq i \leq m-1$, (resp. $\vdash$), (resp. $\dashv$), if the leave points
in the same direction that the $(i-1)^{th}$ (resp. the $(m-1)^{th}$) (resp. $0^{th}$) leave of $\overset{_m}{\treeg}$. Here is an example of how the operations assignement depends on the leaves of here the corolla $\overset{_6}{\treeg}$.
\begin{center}
\includegraphics*[width=4cm]{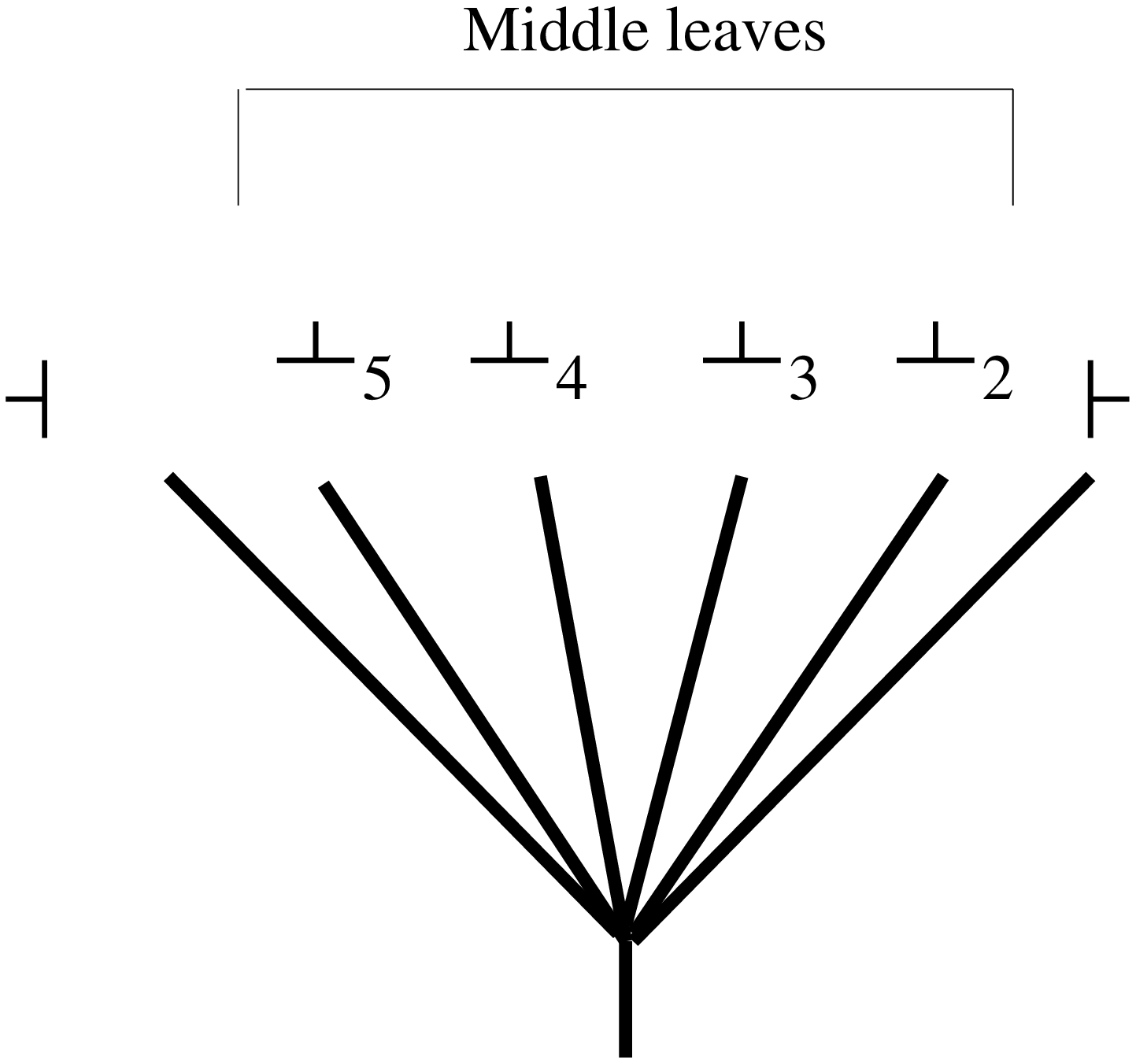}
\end{center}
By convention the last leave of a given $m$-tree remains unassigned. Proceeding that way, a tree from $\overset{_m}{\treeg}_{n+1}$ will give $m$ binary operations, we label 1 to $n$ from left to right.
We have defined for $1\leq j\leq n$, a map,
$$ \circ_j: \overset{_m}{\treeg}_{n+1} \longrightarrow \{\perp_i,\vdash,\dashv \},$$
assigning to each $j$ the corresponding operation
by the process just described. The image will be denoted by $\circ_j^t$. Denote by $\widetilde{\overset{_m}{\treeg}}_n$ the set of trees of $\overset{_m}{\treeg}_n$
whose leaves have been colored by the operations 
$\perp_i,\vdash,\dashv$ just as explained. 
We get a bijection
$tilde: \overset{_m}{\treeg}_n \longrightarrow \widetilde{\overset{_m}{\treeg}}_n$. 
For any $1\leq j \leq n$, define the face map,
$d_j:K\overset{_m}{\treeg}_n \longrightarrow K\overset{_m}{\treeg}_{n-1}$, on $\overset{_m}{\treeg}_n$ first and extended by linearity 
then to be the composite $tilde^{-1} \ \circ del_j \ \circ \ tilde$, where the linear map,
$$del_j: K\widetilde{\overset{_m}{\treeg}}_n \longrightarrow K\widetilde{\overset{_m}{\treeg}}_{n-1},$$
assigns to a tree $t$, the tree $t'$ obtained from $t$ as follows. Localise the leave colored by the operation labelled by $j$. 
Remove the offspring of its father vertex if the $m$ children are leaves, or if all middle leaves are not colored. Otherwise, the result is zero.
Let $T$ be a $m$-tetrahedral algebra over $K$. Define the module of $n$-chains by $C\overset{_m}{\treeg}_n(T):=K[\overset{_m}{\treeg}_n]\otimes T^{\otimes n}$.
Define a linear map $d: C\overset{_m}{\treeg}_n(T) \longrightarrow C\overset{_m}{\treeg}_{n-1}(T)$ by the following formula,
$$ d(t; v_1,\ldots, v_n):=\sum_{j=1}^{n-1} \ (-1)^{j+1}(d_j(t);v_1, \ldots, v_{j-1}, v_j \circ_j^t v_{j+1}, \ldots, v_n):=\sum_{j=1}^{n-1} \ (-1)^{j+1}d_j(t; v_1, \ldots, v_n),  $$
with a slight abuse of notation, and where $t\in \overset{_m}{\treeg}_n$, $v_j \in T$.
\begin{prop}
The face maps $d_l:C\overset{_m}{\treeg}_n(T) \longrightarrow C\overset{_m}{\treeg}_{n-1}(T)$ satisfy the simplicial relations
$d_kd_l=d_{l-1}d_k$ for any $1\leq k < l \leq n-1$.
Therefore $d\circ d=0$ and so $(C\overset{_m}{\treeg}_*(T),d)$ is a chain-complex.
\end{prop}
\Proof
We give a sketch of this idendity for the lowest dimension, that is,
$$d_1d_2=d_{1}d_1:C\overset{_m}{\treeg}_3(T) \rightarrow C\overset{_m}{\treeg}_{2}(T).$$ 
Consider the 5 planar rooted trees on three internal vertices. Remove each offspring of two children by
$m$-children. We obtain 5 $m$-trees of $\overset{_m}{\treeg}_3$ which will give  
dialgebra constraints. Consider now the family of $m$-trees of the form, 
$$((\mid \vee \overset{_m}{\treeg} \vee \ldots \vee \mid)\vee \mid \vee \ldots \vee \mid), \ ((\mid \vee \mid \vee \overset{_m}{\treeg} \vee \ldots \vee \mid)\vee \mid \vee \ldots \vee \mid), \ldots,((\mid \vee \ldots \vee \mid \vee \overset{_m}{\treeg}\vee \mid)\vee \mid \vee \ldots \vee \mid).$$
They will give equalities $(x \perp_{m+1-i}y)\vdash z=0$ for $i \geq 2$.
Apply formally the involution $\dagger$ (compatible with the opposite structure) on these trees to obtain the equalities $x \dashv (y \perp_{i} z)=0$.
Consider now the family of $m$-trees of the form, 
$$(\mid \vee (\overset{_m}{\treeg} \vee \ldots \vee \mid)\vee \mid \vee \ldots \vee \mid), \ (\mid \vee (\mid \vee \overset{_m}{\treeg} \vee \ldots \vee \mid)\vee \mid \vee \ldots \vee \mid), \ldots,(\mid \vee (\mid \vee\ldots \vee \mid \vee \overset{_m}{\treeg})\vee \mid \vee \ldots \vee \mid).$$
We obtain equations $(x\perp_{i} y) \perp_{2} z=0$,
for $i \geq 2$ for the first $m-2$ trees. The last ones gives $x\perp_{2} (y \perp_{2} z)=0$. The $(m-1)^{th}$ one gives $(x\dashv y)\perp_2 z =x\perp_2(y\vdash z)$. Applying the involution $\dagger$ will give the constraints $x\perp_{m-1} (y \perp_{m+1-i} z)=0$, $x\perp_{m-1} (y \perp_{m-1} z)=0$ and $(x\dashv y)\perp_{m-1} z =x\perp_{m-1}(y\vdash z)$.
Similarly the family,
$$(\mid \vee \mid \vee(\overset{_m}{\treeg} \vee \ldots \vee \mid)\vee \mid \vee \ldots \vee \mid), \ (\mid \vee \mid \vee (\mid \vee \overset{_m}{\treeg} \vee \ldots \vee \mid)\vee \mid \vee \ldots \vee \mid), \ldots,$$
$$\ldots (\mid \vee\mid \vee (\mid \vee \ldots \vee \mid \vee \overset{_m}{\treeg})\vee \mid \vee \ldots \vee \mid),$$
will give equations $(x\perp_{i} y) \perp_{3} z=0$,
for $i \geq 3$ for the first $m-3$ trees. The two last one give $x\perp_{3} (y \perp_{2} z)=0$ and
$x\perp_{3} (y \perp_{3} z)=0$. The $(m-2)^{th}$ one gives $(x\dashv y)\perp_3 z =x\perp_3(y\vdash z)$ and so forth.
Observe now that the family,
$$(\mid \vee \ldots \mid \vee \overset{_m}{\treeg} \vee \mid \vee \ldots  \mid \vee  \overset{_m}{\treeg} \vee \mid \vee \ldots \vee \mid),$$
will give the constraints $(x\perp_i y) \perp_j z= x\perp_i (y\perp_j z)$ for $2 \leq i<j \leq m-1$ and also
$(x\perp_i y) \dashv z= x\perp_i (y\dashv z)$ and
$(x\vdash y) \perp_i z= x\vdash (y\perp_i z)$ for $2 \leq i \leq m-1$ as expected.
The general case still splits into two cases. If $j=i+1$, then the proof follows from the low dimension cases and from axioms of $m$-tetrahedral algebras. The case
$j>i+1$ is straightforward.
\eproof

\noindent
We get a chain-complex,
$$C\overset{_m}{\treeg}_*(T):\ \ \ldots \rightarrow K[\overset{_m}{\treeg}_n]\otimes T^{\otimes n} \rightarrow \ldots 
\rightarrow  K[\overset{_m}{\treeg}_3]\otimes T^{\otimes n} \rightarrow K[\overset{_m}{\treeg}_2]\otimes T^{\otimes n} \xrightarrow{\perp_i,\vdash,\dashv} T.$$
This allows us to define the homology of a $m$-tetrahedral algebra $T$ as the homology of the chain-complex $C\overset{_m}{\treeg}_*(T)$, that is $H\overset{_m}{\treeg}_n(T):=H_n(C\overset{_m}{\treeg}_*(T),d)$, $n>0$. The cohomology of $T$ is by definition $H\overset{_m}{\treeg}^n(T):=H^n(C\overset{_m}{\treeg}_*(T),K)$, $n>0$.
For the free $m$-tetrahedral algebra over the $K$-vector space $V$, we get 
$H\overset{_m}{\treeg}_1(m-Tetra.(V))\simeq V$.
We conjecture that for $n>1$,     $H\overset{_m}{\treeg}_n(m-Tetra.(V))=0$, that is the operad $m-Tetra.$ is Koszul.
\section{The Pascal triangle}
We now summarise all the integer sequences we got by gathering them inside the Pascal triangle.
Here is the beginning of this famous triangle. We have chosen the south-east direction (the south-west could be another one) and indicate by an arrow the begining of the coefficients of Poincar\'e series
of the operads involved in this paper and the functors associated with the corresponding categories.
\begin{center}
%\begin{small}
\begin{scriptsize}
$
\begin{array}{ccccccccccccccccc}
& & & & & & & &1 & & & & & & & \\
& & & & & & &1 & &1 & & & & & & &\\
& & & & & & 1& & 2& & 1 & & & &   & &\\
& & & & & 1& & 3 & & 3& &1& & & & &\\
& & & & 1& & 4 & & 6 & &4 & & 1 & & & &\\ 
& & & 1 & & 5 & & 10 & & 10 & & 5 & & 1 & & & \\ 
& & 1& &6& &15 & & 20 & & 15& &6& & 1 & &\\
& \swarrow & & \searrow& &\searrow & & \searrow& &  \searrow & &\searrow  & &\searrow  & &\searrow          
&  \\
As.&\overset{as.}{\leftarrow} & \ldots \hookleftarrow&7-Tetra. & \hookleftarrow& 6-Tetra. & \hookleftarrow& 5-Tetra.&\hookleftarrow & Tetra.&\hookleftarrow &Triang. & \hookleftarrow&Dias. & \hookleftarrow& As.
\end{array}
$
\end{scriptsize}
%\end{small}
\end{center}
Similarly, one could also construct a dual triangle
to the Pascal's one by representing the generalised Catalan numbers for instance as follows. Recall that the operad of associative algebras is self-dual.
\begin{center}
\begin{tiny}
%\begin{scriptsize}
$
\begin{array}{cccccccccccccc}
  & & & & & &1 & & & & & & & \\
  & & & & & 1& & 1& & & & & & \\
  & & & &1 & & 2& & 1 & & & & & \\
  & & & 1& & 3 & 5& 3& &1& & & & \\
  & & 1& & 4 &12 & 14 & 12& 4& & 1 & & & \\ 
 1 & & &5 & 22& 55 & 42& 55 &22 & 5 & & & 1&   \\ 
  \swarrow& &  & \downarrow& \downarrow& \downarrow& \downarrow&\downarrow &  \downarrow &\downarrow & & & & \searrow
  \\
 As. & \overset{+}{\leftarrow} & \ldots & 5-Dend. \hookleftarrow& 4-Dend.\hookleftarrow&3-Dend. \hookleftarrow& 2-Dend.& \hookrightarrow 3-Dend.& \hookrightarrow 4-Dend.&\hookrightarrow 5-Dend. &  & \ldots &\overset{+}{\rightarrow}& As.
\end{array}
$
\end{tiny}
%\end{scriptsize}
\end{center}
\section{Operads whose Poincar\'e series are related to $k$-gonal numbers }
After this long immersion in the realm of triangular shapes and trees, we come back to Section~2.
The free objects associated with the operads $^n\mathcal{P}$ defined in that section seems to be difficult to construct systematically (except for the case $n=2$, see \cite{Loday} and the case $n=3$,
see Section~3).
However, free objects associated with their operadic duals can be computed systematically. Indeed, recall that for $n=3$, the coefficients of the Poincar\'e series of $^3\mathcal{P}:=3-Dend.$ start with $1,3,12,56, \ldots$. Its inverse for the usual composition of functions, that is the Poincar\'e series of the operad $Triang.$, has coefficients starting with $1,3,6,10,15 \ldots$ which are triangular numbers.
Similarly, for $n=4$, the coefficients of the Poincar\'e series of $^4\mathcal{P}$ start with $1,4,23, \ldots$ and seems to be related to noncrossing connected graphs. Its inverse for the usual composition of functions has coefficients starting with $1,4,9,\ldots$ which is the beginning
of square numbers $1,4,9,16,25, \ldots$. For $n=5$,
we could proceed similarly and could conjecture that the inverse series of the Poincar\'e series of $^5\mathcal{P}$, if it exists, would begin by $1,5,12$, which are the beginning
of pentagonal numbers $1,5,12,22, \ldots$, and so forth. Such numbers are generally called $k$-gonal
numbers in the literature (see N.J.A. Sloane Online Encyclopedia of Integers for instance) and the $n^{th}$ $k$-gonal number is, 
$$g_k(n):=n+(k-2)\frac{n(n-1)}{2}.$$
Up to the author's knowedge, serious results on polygonal numbers started with a Fermat's theorem in 1638. He pretended to have proved that any positive integers can always be written as a sum of at most $k$ $k$-gonal numbers. Gauss proved the case of triangular numbers summarising his result by the formula $ E\Upsilon PHKA=\bigtriangleup +\bigtriangleup +\bigtriangleup$ (1796). Euler left
important results on the Fermat's assertion which were used by Lagrange to prove the square case, result also found independently by Jacobi (1772). Finally Cauchy proved the whole assertion (1813).
 
The aim of this part is to construct the operadic dual in the sense of Ginzburg and Kapranov \cite{GK} of
$^k\mathcal{P}$, the associated free objects over a $K$-vector space $V$ and to verify that the coefficients of the Poincar\'e series associated with them are the
$k$-gonal numbers. We start with a definition.
\begin{defi}{}
Fix an integer $k>2$.
A $k$-gonal algebra $G_k$ is a $K$-vector space equipped with $k$ binary operations $\vdash,\dashv,(\perp_i)_{2 \leq i \leq k-1}:G_k^{\otimes 2} \rightarrow G_k$ obeying the following system of quadratic relations for all $2 \leq i,j \leq k-1$ and all $x,y,z \in G_k$, 
\begin{eqnarray}
\begin{cases}
(x\dashv y)\dashv z= x\dashv(y\dashv z),& 
 (x\dashv y)\perp_i z=x \perp_i(y\vdash z),\\
 (x\dashv y)\dashv z= x\dashv(y\vdash z),& 
 (x\vdash y)\perp_i z=x \vdash(y\perp_i z),\\
 (x\vdash y)\dashv z= x\vdash(y\dashv z),& 
 (x\perp_i y)\dashv z=x \perp_i (y\dashv z),\\
(x\dashv y)\vdash z= x\vdash(y\vdash z),& 
(x\perp_i y)\perp_j z=0 =  x \perp_i(y\perp_j z),\\
 (x\vdash y)\vdash z= x\vdash(y\vdash z),& 
 (x\perp_i y)\vdash z=0 = x \dashv(y\perp_i z).
\end{cases}
\end{eqnarray}
The associated category (resp. operad) is denoted by $k$-\textbf{Gonal}, (resp. $k-Gonal$).
\end{defi}
As expected, there are $2(k-2)^2+5(k-1)= 2k^2 -3(k-1)$ relations. 
The functorial diagram between involved categories holds.
\begin{center}
$
\begin{array}{cccc}
\textbf{As.}\rightarrow  &\textbf{Dias.}\rightarrow &\textbf{Triang.}:= 3-\textbf{Gonal} \rightarrow & \textbf{4-Gonal} \rightarrow \ldots \\
  \searrow & \downarrow & \swarrow & \swarrow \\
& \textbf{Leib.}&    &
\end{array}
$
\end{center}
\Rk \textbf{(Opposite $k$-gonal algebra.)}
Let $G_k$ be a $k$-gonal algebra over $K$. Define new operations by:
$$x \vdash' y:= y \dashv x; \ x \dashv' y:= y \vdash x; \ x \perp_i' y:= y \perp_{m+1-i} x,$$
for all $2 \leq i \leq m-1$. Then, the $K$-vector space $G_k$ equipped with these operations is a new
$k$-gonal algebra, denoted by $G_k^{op}$, called the opposite $k$-gonal  algebra. A $k$-gonal is said to be commutative if
 $T^{op}=T$.
\Rk
For any $k$-gonal algebra $G_k$, let $as(G_k)$
be the quotient of $G_k$ by the ideal generated by the elements $x \vdash y - x\dashv y$ and $x \perp_i y$, for all $2 \leq i \leq m-1$ and $x,y \in T$. Observe that $as(G_k)$ is an associative algebra and that the functor $as.(-):$\textsf{$k$-Gonal.}$\rightarrow$ \textsf{As.} is left adjoint to the functor $inc:$ \textsf{As.} $\rightarrow$ \textsf{$k$-Gonal.}.

\noindent
As expected, we get:
\begin{theo}
Fix $k>2$.
The operad \textit{$k$-Gonal} is dual in the sense of \cite{GK} to the operad $^k\mathcal{P}$, that is 
\textit{$k$-Gonal}=$^k\mathcal{P}^!$ and thus $^k\mathcal{P}$=\textit{$k$-Gonal}$^!$.
\end{theo}
To state the next theorem, we introduce for all $1 \leq p \leq k-3$, the linear maps,
$$\psi_1^p: \bigoplus_{n>0} T^{\otimes n} \rightarrow \underset{k-3}{\underbrace{V \oplus \ldots \oplus V}}, \ v_1 \otimes \ldots \otimes v_n \mapsto 0 \oplus \ldots \oplus 0 \oplus \underset{p^{th} \ position}{\underbrace{v_1}} \oplus 0 \oplus \ldots \oplus 0,$$ 
and $\psi_2: \bigoplus_{n>0} T^{\otimes n} \rightarrow T(V)$ defined by $v_1 \otimes \ldots \otimes v_n \mapsto v_2 \otimes\ldots \otimes v_n$ for $n>1$ and by $v \mapsto 1$. 
Denote by $\psi: (K\oplus V) \otimes T(V) \rightarrow K$ and by $\Psi: (K\oplus V) \otimes T(V)\otimes (K\oplus V) \otimes T(V) \rightarrow K$ the canonical projections.
As usual, the free $k$-gonal algebra over a $K$-vector space $V$ is, by definition, the $k$-gonal algebra $k-Gonal.(V)$ equipped with a $K$-linear map
$i:V \hookrightarrow k-Gonal(V)$ such that for any $K$-linear map $f:V \longrightarrow G_k$, where $G_k$ is a $k$-gonal algebra over $K$, there exists a unique $k$-gonal algebra morphism $\phi$ turning the diagram,
\begin{center}
$
\begin{array}{ccc}
&i & \\
V& \hookrightarrow& k-gonal(V) \\
  & & \\
& f \searrow  & \downarrow \phi \\
  & & \\
& & T 
\end{array}
$
\end{center}
commutative.
\begin{theo}
Fix an integer $k>2$.
Let $V$ be a $K$-vector space. Define on 
$$k-Gonal(V):= T(V) \otimes[V \otimes (K \oplus \underset{k-3}{\underbrace{V \oplus \ldots \oplus V}})\otimes T(V)] \otimes T(V),$$
the following binary operations where the $\omega_i \in V^{\otimes p_i}$, $v,v' \in V$ and $w,w'$ belong to either $V$ or $K$.
\begin{footnotesize}
\begin{eqnarray*}
\omega_1\otimes [v\otimes w \otimes \omega_2] \otimes \omega_3 \dashv \omega'_1\otimes [v'\otimes w' \otimes \omega'_2] \otimes \omega'_3 &=&\psi(w',\omega'_2) \omega_1\otimes [v\otimes w \otimes \omega_2]\otimes \omega_3\omega'_1v'\omega'_3,\\
\omega_1\otimes [v\otimes w \otimes \omega_2] \otimes \omega_3 \vdash \omega'_1\otimes [v'\otimes w' \otimes \omega'_2] \otimes \omega'_3 &=&\psi(w,\omega_2)\omega_1v\omega_3\omega'_1\otimes [v'\otimes w' \otimes \omega'_2] \otimes \omega'_3,\\
\omega_1\otimes [v\otimes w \otimes \omega_2] \otimes \omega_3 \perp_2 \omega'_1\otimes [v'\otimes w' \otimes \omega'_2] \otimes \omega'_3 &=&\Psi(w,\omega_2 ,w',\omega'_2)\omega_1\otimes [v
\otimes \omega_3\omega'_1v']\otimes \omega'_3,\\
\omega_1\otimes [v\otimes w \otimes \omega_2] \otimes \omega_3 \perp_i \omega'_1\otimes [v'\otimes w' \otimes \omega'_2] \otimes \omega'_3 &=&\Psi(w,\omega_2 ,w',\omega'_2)\omega_1\otimes [v
\otimes \psi_1^{i-2}(\omega_3\omega'_1v') \otimes \psi_2(\omega_3\omega'_1v')]\otimes \omega'_3,\\
\end{eqnarray*}
\end{footnotesize}
for $3 \leq i \leq k-1$.
Then, $k-Gonal(V)$ is the free $k$-gonal algebra over $V$. Therefore,
the Poincar\'e series of the operad $k-Gonal$ is,
$$ f_{k-gonal}(x)=\sum_{n>0} (-1)^n g_k(n)x^n=\frac{(k-3)x^2-x}{(x+1)^3}.$$
\end{theo}
\Proof
The proof that $k-Gonal(V)$ is a $k$-gonal algebra does not present any difficulties and is left to the reader. The map $i:V \hookrightarrow k-Gonal(V)$ is the composite:
$$ V \simeq K\otimes [V \otimes K \otimes K] \otimes K \hookrightarrow \ \ k-Gonal(V). $$
Let $f:V \longrightarrow G_k$ be a linear map, where $G_k$ is a $k$-gonal algebra. We construct $\phi$ as follows. First on monomials from $k-Gonal(V)$, then extended by $K$-linearity.
Therefore, define $\phi:k-Gonal(V) \longrightarrow G_k$ by:
\begin{footnotesize}
\begin{eqnarray*}
\phi(v_{-r}\otimes \ldots \otimes v_{-1}\otimes &[& v_0\otimes \psi_1^p(w) \otimes u_1\otimes \ldots \otimes u_s ] \otimes v_{1}\otimes \ldots \otimes v_{q}):= \\
& & f(v_{-r})\vdash \ldots \vdash f(v_{-1})\vdash [f(v_0) \perp_{p+2} (f(w) \vdash f(u_1)\vdash \ldots \vdash f(u_s)) ] \dashv f(v_{1})\dashv \ldots \dashv f(v_{q}),\\
\phi(v_{-r}\otimes \ldots \otimes v_{-1}\otimes &[& v_0 \otimes u_1\otimes \ldots \otimes u_s ] \otimes v_{1}\otimes \ldots \otimes v_{q}):= \\
& & f(v_{-r})\vdash \ldots \vdash f(v_{-1})\vdash [f(v_0) \perp_{2} (f(u_1)\vdash \ldots \vdash f(u_s)) ] \dashv f(v_{1})\dashv \ldots \dashv f(v_{q}),
\end{eqnarray*}
\end{footnotesize}
The map $\phi$ is a morphism of $k$-gonal algebra, (use for instance details in the proof of Theorem \ref{freetria}) hence the unicity of such a map since
it has to coincide with $f$ on $V$.
\eproof

\section{Conclusion and open questions}
In the following array, we sum up our results just comming from a simple symmetry on an action of the unit on the most general quadratic relation we could write. 
\begin{center}
\includegraphics*[width=10cm]{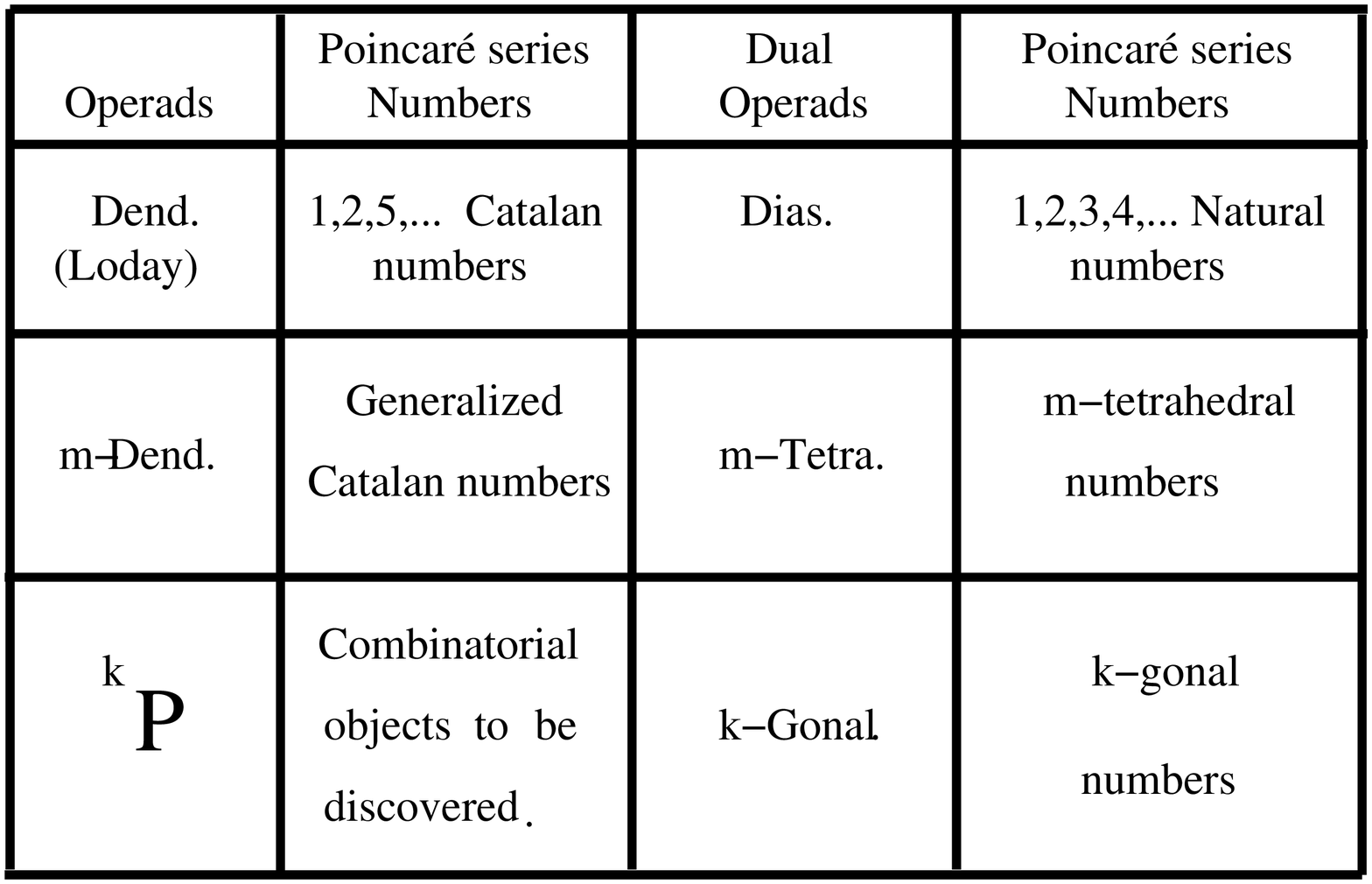}
\end{center}
Some questions come in mind. 
Are these operads considered here Koszul? We conjecture that they are so. Could we propose a theory of additive numbers (like Fermat's theorem) by working only in the operadic setting?
Natural numbers can be manipulated through $Dias.$ and other types of numbers through operads defined
in this paper. On the other hand,
what are the combinatorial objects behind operads $^k\mathcal{P}$, for $k>4$? For instance, 
the inverse of the Poincar\'e series of the operad $5-Tetra.$ is a series whose coefficients
begin with $1,5,38,347,3507,37788,425490, \ldots$. Unfortunately, for the time being no known sequence seems to have this
starting. To end, can we go further on genomics
by using the free involutive 4-dendriform algebra on the generator $\treeCor$?

\bibliographystyle{plain}
\bibliography{These}
L.P. Anciennement rattach\'e \`a l'I.R.M.A.R., Universit\'e de Rennes I et C.N.R.S. U.M.R. 6625 campus de Beaulieu, 35042 Rennes Cedex (France).
\end{document}